\documentclass[9pt,shortpaper,twoside,web]{ieeecolor}
\usepackage{generic}
\usepackage{cite}
\usepackage{amsmath,amssymb,amsfonts}
\usepackage{algorithm}
\usepackage[algo2e]{algorithm2e} 
\usepackage{graphicx}
\usepackage{multicol}
\usepackage{textcomp}
\usepackage{subcaption}
\usepackage{hyperref}

\newtheorem{theorem}{Theorem}[section]
\newtheorem{corollary}{Corollary}[theorem]
\newtheorem{lemma}[theorem]{Lemma}
\newtheorem{proposition}[theorem]{Proposition}

\usepackage{mathtools}
\DeclarePairedDelimiter\floor{\lfloor}{\rfloor}

\DeclareMathOperator*{\argmin}{argmin}

\def\BibTeX{{\rm B\kern-.05em{\sc i\kern-.025em b}\kern-.08em
    T\kern-.1667em\lower.7ex\hbox{E}\kern-.125emX}}
\begin{document}
\title{Parallel and Flexible Dynamic Programming via the\\Randomized Mini-Batch Operator}
\author{Matilde Gargiani, Andrea Martinelli, Max Ruts Martinez and John Lygeros
\thanks{This
work has been supported by the European Research Council (ERC) under the H2020
Advanced Grant no. 787845 (OCAL).
(Corresponding author: Matilde Gargiani.)}
\thanks{M. Gargiani (gmatilde@ethz.ch), A. Martinelli (andremar@ethz.ch) and J. Lygeros (jlygeros@ethz.ch) are with the Automatic Control Laboratory ETH, Physikstrasse 3,
8092 Zurich, Switzerland. }
\thanks{M. R. Martinez (mamartinez@student.ethz.ch) is with ETH Zurich, Rämistrasse 101, 8006 Zurich, Switzerland. }
}

\maketitle

\begin{abstract}
The Bellman operator constitutes the foundation of dynamic programming (DP). An alternative is presented by the Gauss-Seidel operator, whose evaluation, differently from that of the Bellman operator where the states are all processed at once, updates one state at a time, while incorporating into the computation the interim results. The provably better convergence rate of DP methods based on the Gauss-Seidel operator comes at the price of an inherent sequentiality, which prevents the exploitation of modern multi-core systems. In this work we propose a new operator for dynamic programming, namely, the \textit{randomized mini-batch operator}, which aims at realizing the trade-off between the better convergence rate of the methods based on the Gauss-Seidel operator and the parallelization capability offered by the Bellman operator. After the introduction of the new operator, a theoretical analysis for validating its fundamental properties is conducted. Such properties allow one to successfully deploy the new operator in the main dynamic programming schemes, such as value iteration and modified policy iteration. We compare the convergence of the DP algorithm based on the new operator with its earlier counterparts, shedding light on the algorithmic advantages of the new formulation and the impact of the batch-size parameter on the convergence. Finally, an extensive numerical evaluation of the newly introduced operator is conducted. In accordance with the theoretical derivations, the numerical results show the competitive performance of the proposed operator and its superior flexibility, which allows one to adapt the efficiency of its iterations to different structures of MDPs and hardware setups.
\end{abstract}
\section{Introduction}
\label{sec:introduction}
The goal of stochastic optimal control generally is to find a control law that minimizes over a period of time (possibly infinite) a certain measure of the system's performance under uncertainty~\cite{sutton17}. This problem framework is met in numerous applications across different disciplines, from robotics~\cite{bobrow85} to economy~\cite{bertsimas98}, to name a few. Using the concept of value function, R. Bellman introduced in the 1950s what has come to be known as the Bellman equation, a functional equation that allows one to recursively and compactly express in mathematical terms a general stochastic optimal control problem~\cite{bellman52, bellman57}. Dynamic programming (DP) comprises all methods which are based on the solution of the Bellman equation, such as value iteration, policy iteration and their variants. The Bellman equation implicitly defines a functional operator, known as the Bellman operator, whose properties allow one to establish convergence of the various DP methods. Over the years, various extensions of the Bellman operator have been proposed, such as the Gauss-Seidel operator~\cite{bertsekas12} and the relaxed Bellman operator~\cite{martinelli20}. These variants give rise to different DP methods whose convergence depends on the mathematical properties of the underlying operator. Theoretical results and extensive experimentation have confirmed the better convergence rate of DP methods based on the Gauss-Seidel operator~\cite{bertsekas12}. Unfortunately, while the Bellman operator gives rise to DP methods with fully parallelizable iterations, the inherent sequentiality of the Gauss-Seidel operator prevents the exploitation of modern multi-core systems such as GPUs.   

In this paper, we propose a new operator that combines the advantages of the Bellman operator and the Gauss-Seidel operator. In particular, the new operator realises the trade-off between the speed up of parallelization that comes with the Bellman operator and the better convergence rate due to exploitation of interim results that characterizes the Gauss-Seidel operator. In this work the new operator is presented and studied both from a theoretical and practical perspective. 
\noindent The main contributions of this work are the following:
\begin{itemize}
\item Definition of the randomized mini-batch operator, a new operator for DP,
\item A theoretical analysis of its properties as well as the impact of the batch-size on the convergence of the DP method,
\item A publicly available GPU-accelerated Python package that implements value iteration and modified policy iteration based on the new operator.
\end{itemize}
The paper is organized as follows. In Section~\ref{sec:background} we revise the fundamental background of DP and the main DP methods. In Section~\ref{sec:operator} the new operator is introduced and analyzed from a theoretical viewpoint. In particular, first we study its fundamental properties such as shift-invariance, monotonicity and contractivity, and then we study the impact of the batch-size on the convergence of the DP method. Finally, in Section~\ref{sec:experiments} we discuss various numerical results, which corroborate the theoretical results of Section~\ref{sec:operator} and also shed light on the practical advantges of the new operator. 
\section{Background}
\label{sec:background}
\subsection{Problem Framework}
We study stochastic optimal control problems over an infinite horizon and with discounted cost. In particular, our analysis focuses on stationary discrete-time systems with bounded cost per stage and finite state and control spaces. 
More specifically, we consider Markov Decision Processes (MDPs) as a tuple $(\mathcal{S}, \mathcal{U}, P, g, \alpha)$ comprising a finite state space $\mathcal{S}$, a finite control space $\mathcal{U}$, a transition probability function $P: \mathcal{S} \times \mathcal{U} \times \mathcal{S} \rightarrow [0,1]$ that dictates the probability of ending in state $j$ when starting from state $i$ and applying control $u$ in a single time unit, a stage-cost function $g: \mathcal{S} \times \mathcal{U} \rightarrow \mathbb{R}$ that associates a cost to each state-control pair, and a discount factor $\alpha \in (0,1)$. For simplicity we assume that $\mathcal{S} = \left\{1,2,\dots,n \right\}$.
Throughout the paper we use $\mathcal{U}(i)$ to denote the nonempty subset of controls that are associated with state $i$, and $p_{ij}(u)$ to denote the probability of transitioning to state $j$ when applying control $u$ in state $i$. Finally, $g(i,u)$ denotes the stage cost associated with the $(i,u)$-pair. 


Let $\mu:\mathcal{S}\rightarrow\mathcal{U}$, with $\mu(i)\in\mathcal{U}(i)$ be a function that maps states to controls. Such function is generally called \textit{stationary admissible control policy} or, for brevity, \textit{control policy}. We use $\Pi$ to denote the set of all stationary and admissible control policies. In each step $t$ of the decision process, the system is in some state $i$ and a control $u\in\mathcal{U}(i)$ is selected according to a control policy $\mu$ and then applied to the system. The system consequently evolves to a next state $j$ according to the transition probability function. Given an initial state $i_0$, the goal is to find the control policy $\mu$ that minimizes the discounted infinite horizon cost function

\begin{equation}\label{eq: obj_function}
    J_{\mu}(i_0) = \lim_{T\rightarrow \infty}\mathbb{E}_{\mu}\left[\sum_{t=0}^{T-1} \alpha^t g(i_t, \mu(i_t))\right]\,,
\end{equation}
where $(i_0,\mu(i_0), i_1, \mu(i_1), \dots, i_t, \mu(i_t), \dots )$ is the state-control sequence generated by the MDP under policy $\mu$ and starting from state $i_0$, and the expected value is taken with respect to the corresponding probability measure over the space of sequences. The optimal cost function $J^*$ is defined as
\begin{equation}\label{eq: cost-to-go}
    J^*(i) = \min_{\mu\in\Pi} J_{\mu}(i)\,,\quad \forall\,  i \in\mathcal{S}\,.
\end{equation}
A policy $\mu$ is called optimal and denoted with $\mu^*$ if $\mu\in\Pi$ and 
  \begin{equation}
    \label{eq: optimal policy}
     J_{\mu}(i) = J^*(i)\,, \quad \forall\,  i \in \mathcal{S}.
 \end{equation}
Equation~\eqref{eq: obj_function} admits also a \textit{recursive} definition
 
 \begin{equation}\label{eq: bellman equation mu}
     J_{\mu}(i) = g(i,\mu(i)) + \alpha\,  \sum_{j\in\mathcal{S}}p_{ij}(\mu(i))J_{\mu}(j)\,,\quad \forall\,  i \in \mathcal{S} \,,
 \end{equation}
known as the \textit{Bellman equation} for policy $\mu$. Analogously, for the optimal cost $J^*$ we have
  \begin{equation*}
    \begin{aligned}
     J^*(i) &= \min_{\mu\in\Pi}\left[g(i,\mu(i)) + \alpha\,  \sum_{j\in\mathcal{S}}p_{ij}(\mu(i))J^*(j) \right]\,,\quad \forall\,  i \in \mathcal{S}\,.
    \end{aligned}
 \end{equation*}
The minimization over control policies $\mu$ can be written as a minimization over controls $u \in \mathcal{U}(i)$
\begin{equation}
    \label{eq: bellman equation *}
    J^*(i)= \min_{u \in \mathcal{U}(i)}  \left [ g(i, u) + \alpha \sum_{j \in \mathcal{S}} p_{ij}(u) J^*(j)  \right ]\,, \quad \forall\,  i \in \mathcal{S}\,, \\
\end{equation}
which is referred to as the \textit{Bellman equation}. The Bellman equations~\eqref{eq: bellman equation mu} and~\eqref{eq: bellman equation *} play a fundamental role in the development of DP theoretical analysis and the main DP methods~\cite{bertsekas12}. 

\subsection{The Bellman Operators}
Starting from the Bellman equations, it is possible to define two mappings, $T_{\mu}$ and $T$ respectively, known as the \textit{Bellman operators}. These not only provide a compact signature of the problem at hand, but also allow for a convenient shorthand notation for algorithmic description and analysis of DP methods. In particular, given a function $J:\mathcal{S}\rightarrow \mathbb{R}$, we define
\begin{equation}\label{eq: Bellman operator mu}
   \left( T_{\mu}J \right)(i) =  g(i,\mu(i)) + \alpha\sum_{j\in\mathcal{S}} p_{ij}(\mu(i))J(j)\,,\quad \forall\,  i \in \mathcal{S} \,,
\end{equation}
and
\begin{equation}\label{eq: Bellman operator *}
   (TJ)(i) = \min_{u\in \mathcal{U}(i)}\left[ g(i,u) + \alpha\sum_{j\in\mathcal{S}} p_{ij}(u)J(j) \right]\,,\quad \forall\,  i \in \mathcal{S}\,.
\end{equation}
The Bellman operators~\eqref{eq: Bellman operator mu} and~\eqref{eq: Bellman operator *} allow one to rewrite the Bellman Equations~\eqref{eq: bellman equation mu} and~\eqref{eq: bellman equation *} as the fixed point equations $J_{\mu} = T_{\mu}J_{\mu}$ and $J^* = TJ^*$, respectively.
It can be shown that the Bellman operators are contractive, so, thanks to the Banach Theorem~\cite{rockafellar76}, they admit unique fixed points $J_{\mu}$ and $J^*$, respectively. Moreover, the corresponding Picard-Banach iterations converge to the fixed points from any initial condition $J$

\begin{minipage}{0.2\textwidth}
  \begin{equation}
    	\label{eq: convergence of T_mu}
    	\lim_{k\rightarrow \infty} T_\mu^k J = J_\mu\,,
  	\end{equation}
\end{minipage}\quad\quad\quad
\begin{minipage}{0.2\textwidth}
    \begin{equation}
    	\label{eq: convergence of T}
    	\lim_{k\rightarrow \infty} T^k J = J^*\,.
  	\end{equation}
\end{minipage}

\subsection{The Main DP Methods}
Equations~\eqref{eq: convergence of T_mu} and~\eqref{eq: convergence of T} allow the Bellman operators to be deployed in the main computational DP schemes, such as value iteration (VI) and modified policy iteration (MPI). Value iteration, also known as the method of \textit{successive approximations}, is based on Equation~\eqref{eq: convergence of T}. The method consists in starting from an arbitrary finite cost $J$ and then repeatedly applying the operator $T$ to generate a sequence of refined estimates that in the limit converges to the optimal cost $J^*$. To retrieve the optimal policy $\mu^*$, one can refer to Equation \eqref{eq: optimal policy}, which states that an optimal policy is greedy with respect to the optimal cost function $J^*$. 

An alternative method for directly obtaining the optimal policy is provided by policy iteration (PI). 
In contrast to VI, the PI algorithm converges to the optimal policy $\mu^*$ and the optimal cost $J^*$ in a finite number of steps, since the policy is improved at each iteration and since, by the finiteness of $\mathcal{S}$ and $\mathcal{U}$, there only exists a finite number of stationary policies $\mu$. Nevertheless, the PI algorithm has the disadvantage that, in order to compute the exact cost function of a policy, a system of linear equations has to be solved. The dimension of such system is equal to the number of states and therefore, for large state-spaces, may lead to a very expensive computation. An alternative scheme is provided by the modified policy iteration algorithm described in Algorithm~\ref{alg: MPI}.
\begin{algorithm}[H]
\caption{Modified Policy Iteration}
\label{alg: MPI}
\SetAlgoLined
 \textit{Initialization:}
 
 $\mu \in \Pi,\,K\in \mathbb{N}$
 
\Repeat{policy has converged}{
     
     \textit{Policy evaluation:}
    
     $J(i) \in \mathbb{R}, \quad\forall\,  i \in \mathcal{S}$
     
     \For{$k=1,2,\dots,K$}{
     
    
    $J(i) \leftarrow T_\mu J(i), \quad \forall\,  i \in \mathcal{S}$}
    
     \textit{Policy improvement:}
    
     $\mu(i)\!\! \leftarrow \!\!\argmin_{u \in \mathcal{U}(i)}\!\!  \left[ g(i, u) + \alpha \sum_{j \in \mathcal{S}} p_{ij}(u) J(j)\right ], \forall\,  i \in \mathcal{S} $
}
 
\end{algorithm}
\noindent The main difference with respect to standard PI is that the policy evaluation step is carried out approximately at every iteration by applying $K$ iterations of the VI method~\eqref{eq: convergence of T_mu} with the $T_{\mu}$ operator, starting from an arbitrary estimate $J$ of the cost associated to the current policy $\mu$. This version of the PI algorithm is therefore more appealing than the standard PI for systems with large state-spaces, since it does not require the exact solution of a potentially large system of linear equations. This inexact variant of PI still allows one to recover the optimal policy in a finite number of iterations. If warm-starting is deployed in the policy evaluation step, then the method also asymptotically converges to the optimal cost $J^*$. 
See Chapter 2 in~\cite{bertsekas12} for a more thorough description of VI and MPI.

\subsection{The Gauss-Seidel Operators}

Note that the computation of the Bellman operators can be carried out in parallel for each state, as they do not make use of interim results. An alternative approach is presented by the \textit{Gauss-Seidel operators}, whose computation consists in updating one state at a time, while incorporating the interim results. This alternative operator is inspired by the Gauss-Seidel method for solving linear and nonlinear systems of equations~\cite{kelley95}. In particular, the Gauss-Seidel operator $F$ is defined as
\begin{equation*}
    FJ(1) = \min_{u \in \mathcal{U}(1)} \!\left [ g(1,u) + \alpha \sum_{j=1}^n p_{1j}(u)J(j) \right ],
\end{equation*}
and for $i=2,\dots,n$
\begin{equation*}
    FJ(i)\! = \!\!\min_{u \in \mathcal{U}(i)} \left [ g(i,u) + \alpha\! \sum_{j=1}^{i-1} p_{ij}(u)FJ(j) + \alpha \!\sum_{j=i}^n p_{ij}(u)J(j) \right ].
\end{equation*}
To simplify the notation, one can also introduce the set of the states that have already been updated when processing state $i$
\begin{equation*}
     \mathcal{M}_F(i) = \left\{1,2,\dots,i-1 \right\},\quad \forall i \in \mathcal{S}\,,
\end{equation*}
leading to
\begin{equation*}
    FJ(i)\! =\!\! \min_{u \in \mathcal{U}(i)} \!\!\left [ g(i,u) + \alpha\!\!\!\!\!\!\!\!\!\!\! \sum_{j \in \mathcal{S}\setminus \mathcal{M}_F(i)}\!\!\!\!\!\!\!\!\! p_{ij}(u)J(j) + \alpha\! \!\!\!\!\!\!\! \sum_{j \in \mathcal{M}_F(i)} \!\!\!\!\!\!\! p_{ij}(u)FJ(j) \right ]\,.
\end{equation*}
The Gauss-Seidel operator $F_{\mu}$ for policy evaluation can be easily deduced and therefore it is omitted in the interest of space.

Like the Bellman operators, the Gauss-Seidel operators are also $\alpha$-contractive leading to equations analogous to \eqref{eq: convergence of T_mu} and \eqref{eq: convergence of T} for the Gauss-Seidel case~\cite{bertsekas12}. By deploying the $F$ operator in place of the $T$ operator in the VI algorithm, the Gauss-Seidel version of VI (GS-VI) is obtained.
Similarly, by utilizing the $F_\mu$ operator in place of the $T_\mu$ operator in the policy evaluation step of Algorithm~\ref{alg: MPI}, the Gauss-Seidel version of MPI (GS-MPI) is obtained. Note that the computation of $FJ$ and $TJ$ has the same complexity while well-enstablished theoretical results and extensive experimentation indicate the better convergence rate of GS-VI and GS-MPI~\cite{bertsekas12}. The better convergence rate though comes with an inherent sequentiality. In fact, in contrast to the Bellman operators, the computation of the Gauss-Seidel operators can not be parallelized, since the computation of the map for one state requires the computation of the map for all the previous states preventing the exploitation of modern multi-core systems. As a consequence, despite their better convergence rate, GS-VI and GS-MPI might require more computation time with respect to their Bellman counterparts on a parallel system. 

\section{The Randomized Mini-Batch Gauss-Seidel Operator}\label{sec:operator}
With the aim of realizing the trade-off between better convergence rate and parallelization capability, and taking inspiration from the block Gauss-Seidel method~\cite{bertsekas99,beck13} and the mini-batch stochastic optimization methods~\cite{richtarik2011,bottou18}, we propose a new set of operators, the \textit{randomized mini-batch operators}, which, for brevity, we also call \textit{mini-batch operators}. The new operators update the states in batches, while incorporating into the computation the interim results of states belonging to previously updated batches. As a result, the computation within  batches is fully parallelizable and this allows the user to select the batch-size that best exploits the parallelization capabilities of the hardware at hand. At the same time, the convergence rate benefits from the use of interim results of states belonging to previously updated batches. In this section, the mini-batch operators are first defined as a generalization of the Gauss-Seidel and Bellman operators. Then their deployment in the main DP methods is briefly discussed. Finally, we conduct a theoretical analysis of their fundamental properties and compare the convergence of the DP method based on the mini-batch operator with that of its counterparts.

\subsection{Definition}
Without loss of generality, we assume that the states are processed in ascending order. In practice, one can first assign random indexes to the states and then process them in an ascending order; indeed, re-randomisation can be performed for each iteration for VI and MPI (see Section~\ref{sec:experiments}). Let $m\in\mathbb{N}$, $1\leq m\leq n$ be the batch-size. Analogously to the Gauss-Seidel operators, we define the function $\mathcal{M}_{m} :\mathcal{S}\rightarrow 2^{\mathcal{S}}$ as follows
\begin{equation}
\label{eq: M(i) definition}
\mathcal{M}_m(i) = \bigcup_{j=1}^{m \cdot \floor{\frac{i-1}{m}}} \left \{ j \right \}, \quad \forall\,  i \in \mathcal{S}\,,
\end{equation} 
where $\mathcal{M}_{m}(i)$ is a function parameterized by the batch-size $m$ that maps state $i$ to the set of states that have been already updated when processing state $i$. This function allows one to define a general operator parameterized by the batch-size $m$, which facilitates the theoretical analysis and provides a unified framework for comparing different operators.
In particular, given a function $J:\mathcal{S}\rightarrow \mathbb{R}$, the mini-batch operators $B_{m}$ and $B_{\mu, m}$ are defined as follows 
\begin{equation}
\begin{split}
\label{eq: B definition}
    B_m J(i) = \min_{u \in \mathcal{U}(i)} \Bigg [ g(i,u) &+ \alpha \!\!\!\!\!\!\!\!\!\! \sum_{j \in \mathcal{S}\setminus \mathcal{M}_m(i)} \!\!\!\!\!\!\!\! p_{ij}(u)J(j) \\ &  +\, \alpha \!\!\!\!\!\!\!\! \sum_{j \in \mathcal{M}_m (i)} \!\!\!\!\!\! p_{ij}(u)B_m J(j) \Bigg]\,,
\end{split}
\end{equation}

\begin{equation}
\begin{split}
\label{eq: Bmu definition}
B_{\mu,m} J(i) =  g(i,\mu(i)) &+ \alpha \!\!\!\!\!\!\!\!\!\! \sum_{j \in \mathcal{S}\setminus \mathcal{M}_m(i)} \!\!\!\!\!\!\!\! p_{ij}(\mu(i))J(j) \\ &+ \alpha \!\!\!\!\!\!\!\! \sum_{j \in \mathcal{M}_m (i)} \!\!\!\!\!\! p_{ij}(\mu(i))B_{\mu,m} J(j)\,.
\end{split}
\end{equation}

Notice that, if $m=n$, then $\mathcal{M}_n(i) = \left\{\emptyset\right\}$ for all $i\in\mathcal{S}$ and the Bellman operators are recovered. If $m=1$, then $\mathcal{M}_1(i)$ is equivalent to $\mathcal{M}_F(i)$ for all $i\in\mathcal{S}$ and the Gauss-Seidel operators are recovered. To simplify the notation, unless needed, we neglect the dependence on the batch-size $m$ and use the more compact notation $B$, $B_{\mu}$ and $\mathcal{M}(i)$.

\subsection{Theoretical Analysis}
By deploying the mini-batch operators in place of the Bellman operators in the VI and MPI schemes we obtain the mini-batch versions of VI (MB-VI) and MPI (MB-MPI). The theoretical analysis of the properties of these algorithms is divided into two parts: first, in Section~\ref{fp and c}, we analyze the fundamental properties of the mini-batch operators, \textit{i.e.}, monotonicity (Lemma~\ref{lmm: monotonicity}), shift-invariance (Lemma~\ref{lmm: shift invariance}) and contractivity (Proposition~\ref{prp: alpha-contractivity}). These properties are then used to characterize the fixed points of the operators (Lemma~\ref{lmm: existence and uniqueness}) as well as the convergence of the DP method based on the mini-batch operators (Theorem~\ref{prp: convergence}). In the second part (Section~\ref{sec: comparative analysis}), we provide a comparative analysis of the mini-batch operators which sheds light on the the impact of the batch-size on the convergence rate (Theorem~\ref{thm: comparison} and Corollary~\ref{cr: comparison}). In the interest of space, the proofs are carried out only for $B$, but similar derivations are applicable also for $B_{\mu}$.
\subsubsection{Fundamental Properties \& Convergence}\label{fp and c}
The next two lemmas and proposition characterize the fundamental properties of the mini-batch operators, namely, monotonicity, shift-invariance and $\alpha$-contractivity in infinity norm. 
\begin{lemma}[monotonicity]
\label{lmm: monotonicity}
For any two functions $J:\mathcal{S}\rightarrow\mathbb{R}$ and  $J':\mathcal{S}\rightarrow\mathbb{R}$ and for any stationary policy $\mu$, if $J\leqslant J'$, then $B^kJ \leqslant B^kJ'$ and $B_{\mu}^kJ \leqslant B_{\mu}^kJ'$ for $k=1,2,\dots$.

\begin{proof}
We use induction twice, first on the state counter $i$ and then on the iteration counter $k$. We first set $k=1$ and apply induction over the state counter $i$; to simplify notation we drop the dependency on $k$.
As base-case, we consider $i=1$ and apply the $B$ operator to $J$. Since $\mathcal{M}(1)=\{\emptyset\}$, Equation~\eqref{eq: B definition} implies   
\begin{equation*}
\begin{split}
     (BJ)(1) & = \min_{u \in \mathcal{U}(1)}  \left [ g(1, u) + \alpha \sum_{j \in \mathcal{S}} p_{1j}(u) J(j)  \right ]\\
     & \leq \!\min_{u \in \mathcal{U}(1)}\!  \left [ g(1, u) + \alpha \sum_{j \in \mathcal{S}} p_{1j}(u) J'(j)  \right ] \\
     & = (BJ')(1)\,,\\
     \end{split}
\end{equation*}    
where the inequality follows from the assumption that $J\leq J'$.
We then consider an arbitrary state index $i>1$ and we assume that $\left( BJ\right)(j)\leq\left( BJ'\right)(j)$ for $j=1,2,\dots,i-1$.
Since, without loss of generality, all states are assumed to be processed in ascending order and since all states in $\mathcal{M}(i)$ have already been processed by definition, then $\mathcal{M}(i) \subseteq \left\{1,2,\dots,i-1 \right\}$. Therefore the induction assumption implies
\begin{equation}
\label{eq: BJ leq BJ'}
\left( BJ\right)(j)\leq\left( BJ'\right)(j),\quad \forall\, \,j\in\mathcal{M}(i)\,. 
\end{equation} 
By applying the $B$ operator to $J$ for state $i$ and considering Equation~\eqref{eq: BJ leq BJ'} together with the fact that $J\leq J'$, we obtain
\begin{equation*}
\begin{split}
     (BJ)&(i)\\
     &=\!\! \min_{u \in \mathcal{U}(i)}\! \left [ g(i, u) + \alpha \!\!\!\!\!\!\!\!\!\sum_{j \in \mathcal{S} \setminus \mathcal{M}(i)} \!\!\!\!\!\!\!\! p_{ij}(u) J(j) + \alpha \!\!\!\!\!\! \sum_{j \in \mathcal{M}(i)}\!\!\!\!\! p_{ij}(u) (B J)(j) \right] \\
     & \leq \!\!\! \min_{u \in \mathcal{U}(i)} \!\! \left [\! g(i, u) + \alpha\!\!\!\!\!\!\!\!\! \sum_{j \in  \mathcal{S} \setminus \mathcal{M}(i)} \!\!\!\!\!\!\!\! p_{ij}(u) J'(j) + \alpha\!\!\!\!\!\! \sum_{j \in \mathcal{M}(i)}\!\!\!\!\! p_{ij}(u) (B J')(j) \!\right] \\
       & = (BJ')(i),
\end{split}
\end{equation*}
which is a valid inductive step over $i$.
We can therefore conclude that 
\begin{equation}
\label{eq: base-case on k}
\left(BJ\right)(i) \leq (BJ')(i),\quad i = 1,\dots, n\,. 
\end{equation}

We now apply induction over the iteration counter $k$. The base-case is given by Equation~\eqref{eq: base-case on k}.
We assume that the following inequalities hold for an arbitrary integer $k\geq 1$ 
\begin{equation}
\label{eq: BkJ leq BkJ'}
( B^kJ)(i) \leq ( B^kJ')(i),\quad i=1,\dots,n \,.
\end{equation}
Without loss of generality, we introduce the costs $\tilde{J}$ and $\tilde{J}'$ such that $\tilde{J}=B^kJ$ and $\tilde{J}'=B^kJ'$. Consequently, Inequality~\eqref{eq: BkJ leq BkJ'} can be equivalently rewritten as $\tilde{J} \leq \tilde{J}'$. By applying the $B$ operator on both sides of Inequality~\eqref{eq: BkJ leq BkJ'} and making use of the monotonicity property of $B$ (Lemma~\ref{lmm: monotonicity}), we obtain that for $i=1,\dots,n$ 
\begin{equation*}
\begin{aligned}
 (B^{k+1}J)(i) &= B (B^k J)(i) \\
 &= (B\tilde{J})(i) \\
 &\leq (B\tilde{J'})(i) \\
 &= B(B^kJ')(i)\\
 & = (B^{k+1}J')(i)\,,
\end{aligned}
\end{equation*}
which is a valid induction step over $k$.
We can therefore conclude that $(B^kJ)(i) \leqslant (B^kJ')(i)$ for $i = 1,\dots,n$ and $k = 1,2,\dots$.
\end{proof}
\end{lemma}
\begin{lemma}[shift-invariance]
\label{lmm: shift invariance}
For any function $J: \mathcal{S} \rightarrow \mathbb{R}$, stationary policy $\mu$ and $r\in\mathbb{R}$, then $(B^k(J + re))(i) = (B^kJ)(i) + \alpha^k r$ and $(B_{\mu}^k(J + re))(i) = (B_{\mu}^kJ)(i) + \alpha^k r$ for $i = 1,\dots,n$ and $k=1,2,\dots$, where $e$ is the unit function that takes value 1 identically on $\mathcal{S}$.

\begin{proof}
Similarly to Lemma~\ref{lmm: monotonicity}, we use induction twice, first on the state counter $i$ and then on the iteration counter $k$. We start by setting $k=1$ and applying induction over the state counter $i$; to simplify notation we drop the dependency on $k$.
As base-case we set $i=1$ and then we apply the $B$ operator to $J + r e$. Since $\mathcal{M}(1)=\{\emptyset\}$, Equation~\eqref{eq: B definition} implies    
\begin{equation*}
\begin{split}
(B(J + &re))(1) \\ & = \min_{u \in \mathcal{U}(1)}  \left [ g(1, u) + \alpha \sum_{j \in \mathcal{S}} p_{1j}(u) (J + re)(j) \right ] \\
   & = \min_{u \in \mathcal{U}(1)} \left [ g(1, u) + \alpha \sum_{j \in \mathcal{S}} p_{1j}(u) J(j) +  \alpha \sum_{j \in \mathcal{S}} p_{1j}(u) r \right ] \\
   & = \min_{u \in \mathcal{U}(1)} \left [ g(1, u) + \alpha \sum_{j \in \mathcal{S}} p_{1j}(u) J(j) \right ] + \alpha r \\
   & = (BJ)(1) + \alpha r\,.
\end{split}
\end{equation*}
Then we consider an arbitrary state index $i>1$ and assume that $(B(J + r e))(j) = (BJ)(j) + \alpha r$ for $j=1,2,\dots,i-1$. 
Since, without loss of generality, all states are assumed to be processed in ascending order and since all states in $\mathcal{M}(i)$ have already been processed by definition, then $\mathcal{M}(i) \subseteq \left\{1,2,\dots,i-1 \right\}$. Therefore the induction assumption implies
\begin{equation}
\label{assumption shift invariance induction for b}
    (B(J + re))(j) = (BJ)(j) + \alpha r, \quad \forall\,  j \in \mathcal{M}(i)\,.
\end{equation}
By applying the $B$ operator to $J + re$ for state $i$ and considering Equation~\eqref{assumption shift invariance induction for b}, we obtain
\begin{equation*}
\begin{split}
    &(B(J + re))(i) \\ &=  \min_{u \in \mathcal{U}(i)} \!\!\Bigg[ g(i, u) + \alpha \!\! \mkern-22mu \sum_{j \in  \mathcal{S} \setminus \mathcal{M}(i)} \mkern-15mu\!\! p_{ij}(u) (J + re)(j) \\ 
  &   \phantom{ \min_{u \in \mathcal{U}(i)} \!\!\Bigg[ g(i, u) }\quad\quad\quad\quad\quad\,\,\,\,\,\,\,\,\,\,\, + \alpha\!\! \mkern-13mu \sum_{j \in \mathcal{M}(i)} \mkern-9mu\!\! p_{ij}(u) (B(J + re))(j) \Bigg] \\
   & \mkern-5mu\stackrel{(\ref{assumption shift invariance induction for b})}{=} \min_{u \in \mathcal{U}(i)}  \!\!\Bigg[ g(i, u) + \alpha\!\! \mkern-22mu \sum_{j \in \mathcal{S} \setminus \mathcal{M}(i)} \mkern-15mu \!\!p_{ij}(u) J(j) + \alpha \!\!\mkern-13mu \sum_{j \in  \mathcal{M}(i)} \mkern-8mu\!\! p_{ij}(u) (B J)(j)  \\
   &  \phantom{ \min_{u \in \mathcal{U}(i)} \!\!\Bigg[ g(i, u)  + \alpha\!\! \mkern-22mu \sum_{j \in \mathcal{S} \setminus \mathcal{M}(i)} \mkern-15mu \!\!p_{ij}(u) J(j) }\quad\quad\,\,\,\,\,\,\,\,\,\,\,\,\,\, + \alpha \mkern-5mu \sum_{j \in \mathcal{S}} p_{ij}(u) r\Bigg]  \\
   & = \min_{u \in \mathcal{U}(i)}\!\! \left [ g(i, u) + \alpha \!\!\mkern-22mu \sum_{j \in \mathcal{S} \setminus \mathcal{M}(i)} \mkern-15mu\!\! \! p_{ij}(u) J(j) + \alpha \!\!\mkern-13mu \sum_{j \in  \mathcal{M}(i)} \mkern-10mu\! p_{ij}(u) (B J)(j) \right ] \!+\! \alpha r \\
   & = (BJ)(i) + \alpha r,
\end{split}
\end{equation*}
which is a valid inductive step over $i$.
We can therefore conclude that 
\begin{equation}
  \label{eq: shift invariance result of induction for b}
   (B(J + re))(i) = (BJ)(i) + \alpha r\,, \quad i = 1,\dots, n\,.
\end{equation}

We now apply induction over the iteration counter $k$. The base-case is given by Equation~\eqref{eq: shift invariance result of induction for b}.
We then assume that the following inequalities hold for an arbitrary integer $k\geq 1$ 
\begin{equation}
\label{eq: induction assumption k shift}
(B^k(J + r e))(i) = (B^k J)(i) + \alpha^k r\,,\quad i=1,\dots,n\,.
\end{equation}
Analogously to Lemma~\ref{lmm: monotonicity}, we introduce $\tilde{r}$ and $\tilde{J}$ where $\tilde{r} = \alpha^k r$ and $\tilde{J} = B^k J$. By applying the $B$ operator on the left-hand side of Equation~\eqref{eq: induction assumption k shift} and making use of Equation~\eqref{eq: shift invariance result of induction for b}, we obtain
\begin{equation*}
\begin{aligned}
(B^{k+1}(J + r e))(i) &= (B(B^k(J + r e) ))(i)\\
& \mkern-5mu \stackrel{(\ref{eq: induction assumption k shift})}{=} (B(B^k J + \alpha^k r e))(i)\\
& = (B(\tilde{J} + \tilde{r}e))(i)\\
& \mkern-5mu\stackrel{(\ref{eq: shift invariance result of induction for b})}{=} (B\tilde{J})(i) + \alpha \tilde{r}\\
&= (B(B^k J))(i) + \alpha\alpha^k r\\
& = (B^{k+1}J)(i) + \alpha^{k+1}r\,,
\end{aligned}
\end{equation*}
which is a valid induction step over $k$.
We can therefore conclude that $(B^k(J + re))(i) = (B^kJ)(i) + \alpha^k r$ for $i = 1,\dots,n$ and $k = 1,2,\dots$.
\end{proof}
\end{lemma}
\begin{proposition}[$\alpha$-contractivity in infinity norm]
\label{prp: alpha-contractivity}
For any two functions $J:\mathcal{S}\rightarrow \mathbb{R}$ and $J':\mathcal{S}\rightarrow \mathbb{R}$, then for $k=1,2,\dots$
\begin{equation*}
\max_{i \in \mathcal{S}} |(B^kJ)(i)-(B^kJ')(i)| \leqslant \alpha^k\!\! \max_{i \in \mathcal{S}} |J(i)-J'(i)|\,,
\end{equation*}
\begin{equation*}
\max_{i \in \mathcal{S}} |(B_\mu^kJ)(i)-(B_\mu^kJ')(i)| \leqslant \alpha^k\!\! \max_{i \in \mathcal{S}} |J(i)-J'(i)|\,.
\end{equation*}

\begin{proof}
Let $c = \max_{{i \in \mathcal{S}}} \vert J(i) - J'(i) \vert$. Then
\begin{equation}
\label{eq: contraction inequality}
J(i) - c \leq J'(i)\leq J(i) + c\,,\quad i=1,\dots,n\,.
\end{equation}
By applying $B^k$ to each inequality in Equation~\eqref{eq: contraction inequality} and using Lemma~\ref{lmm: monotonicity} and Lemma~\ref{lmm: shift invariance}, we obtain that for $i=1,\dots,n$ 
\begin{equation}
\label{eq: set of ineq contractivity}
\begin{aligned}
(B^k(J-c e))(i) &= (B^kJ)(i) - \alpha^k c\\
&\leq (B^k J')(i)\\
&\leq (B^k(J + c e))(i)\\
&= (B^k J')(i) + \alpha^k c\,.
\end{aligned}
\end{equation}
Equation~\eqref{eq: set of ineq contractivity} reduces to
\begin{equation*}
    (B^kJ)(i) - \alpha^k c \leqslant (B^kJ')(i) \leqslant  (B^kJ)(i) + \alpha^k c, \,\, i =1,\dots,n\,,
\end{equation*}
which can be equivalently reformulated as
\begin{equation}
\label{eq: abs-value ineq contraction}
    |(B^kJ)(i)-(B^kJ')(i)| \leqslant \alpha^k c, \quad i =1,\dots,n\,.
\end{equation}
Finally, Equation~\eqref{eq: abs-value ineq contraction} trivially implies that
\begin{equation*}
    \max_{i \in \mathcal{S}} |(B^kJ)(i)-(B^kJ')(i)| \leqslant \alpha^k \max_{i \in \mathcal{S}} |J(i)-J'(i)|\,,
\end{equation*}
which concludes the proof.
\end{proof}
\end{proposition}
The following lemma characterizes the unique fixed points of the mini-batch operators. In particular, the $B$ operator has $J^*$ as unique fixed point, while the unique fixed point of the $B_{\mu}$ operator is $J_{\mu}$.
\begin{lemma}[existence \& uniqueness of fixed points]\label{lmm: existence and uniqueness}
The optimal cost function $J^*$ and the cost $J_{\mu}$ associated with the policy $\mu$ are the unique fixed points of the operators $B$ and $B_{\mu}$, respectively.  

\begin{proof}
To prove that $J^*$ is a fixed point of the $B$ operator, we apply induction over the state counter $i$. As base-case we set $i=1$ and apply the $B$ operator to $J^*$ for state $i=1$. Since $\mathcal{M}(1)=\{\emptyset\}$, we obtain
\begin{equation*}
\begin{aligned}
(BJ^*)(1) &=\min_{u \in \mathcal{U}(1)}  \left [ g(1, u) + \alpha \sum_{j \in \mathcal{S}} p_{1j}(u) J^*(j) \right]\\
        &\stackrel{(\ref{eq: bellman equation *})}{=} J^*(1)\,.
\end{aligned}
\end{equation*}
We then consider an arbitrary state index $i> 1$ and assume that $(BJ^*)(j) = J^*(j)$ for $j=1,2,\dots,i-1$.
Since, without loss of generality, all states are assumed to be processed in ascending order and since all states in $\mathcal{M}(i)$ have already been processed by definition, then $\mathcal{M}(i) \subseteq \left\{1,2,\dots,i-1 \right\}$. Therefore the induction assumption implies
\begin{equation}
\label{eq: assumption induction contraction}
    (B J^*)(j) = J^*(j)\,, \quad \forall\,  j \in \mathcal{M}(i)\,.
\end{equation}
By applying the $B$ operator to $J^*$ for state $i$ and taking into account Equation~\eqref{eq: assumption induction contraction}, we obtain
\begin{equation*}
\begin{aligned}
\!(BJ^*)(i) & \!=\!\!\!\! \min_{u \in \mathcal{U}(i)}  \!\!\!\left[ \!g(i, u) \!+\! \alpha\!\!\! \mkern-20mu\sum_{j \in \mathcal{S} \setminus \mathcal{M}(i)}\mkern-18mu\!\! p_{ij}(u) J^*(j)\! +\! \alpha\!\! \mkern-13mu\sum_{j \in \mathcal{M}(i)} \mkern-8mu\!\! p_{ij}(u) (B J^*)(j)\! \right] \\
        & \mkern-5mu\stackrel{(\ref{eq: assumption induction contraction})}{=} \min_{u \in \mathcal{U}(i)}  \left [ g(i, u) + \alpha \sum_{j \in \mathcal{S}} p_{ij}(u) J^*(j) \right]\\
        &\mkern-3mu\stackrel{(\ref{eq: bellman equation *})}{=} J^*(i)\,,
\end{aligned}
\end{equation*}
which is a valid inductive step over $i$.
We can conclude that $(BJ^*)(i) = J^*(i)$ for $i=1,\dots,n$.
Uniqueness follows directly from the Banach Theorem thanks to $\alpha$-contractivity.
\end{proof}
\end{lemma}
The following proposition shows that the DP method based on the mini-batch operators converges to their unique fixed points, \textit{i.e.}, $J^*$ and $J_{\mu}$, starting from any arbitrary bounded function $J$.
\begin{proposition}[convergence of the mini-batch DP method]
\label{prp: convergence}
For any function $J:\mathcal{S}\rightarrow \mathbb{R}$ and stationary policy $\mu$, then $J^* = \lim_{k \to \infty} B^kJ$ and $J_\mu = \lim_{k \to \infty} B_\mu^kJ$ for $k=1,2,\dots$.

\begin{proof}
Convergence follows directly from Proposition~\ref{prp: alpha-contractivity} and Lemma~\ref{lmm: existence and uniqueness}.
In particular, we start from Equation~\eqref{eq: abs-value ineq contraction} and substitute $J^*$ in place of $J'$ as follows
\begin{equation*}
\begin{aligned}
\vert (B^k J)(i) - (B^k J^*)(i) \vert &= \vert (B^k J)(i) - J^*(i) \vert  \\
&\leq \!\alpha^k \!\max_{j\in\mathcal{S}}\vert J(j) - J^*(j) \vert\,,\,\, i=1,\dots,n\,.
\end{aligned}
\end{equation*}
Then, taking the limit as $k\rightarrow \infty$, we obtain the upper bound
\begin{equation*}
\begin{aligned}
\lim_{k\rightarrow\infty}  \vert (B^k J)(i) - J^*(i) \vert &\leq \lim_{k \rightarrow \infty} \alpha^k \vert J(i) - J^*(i) \vert = 0 \,,
\end{aligned}
\end{equation*}
which implies that $\lim_{k\rightarrow\infty} (B^k J)(i) = J^*(i)$ for $i=1,\dots,n$.
\end{proof}
\end{proposition}

\subsubsection{Comparative Analysis}
\label{sec: comparative analysis}
The following theorem and corollary characterize the impact of the batch-size on the convergence rate of the DP method based on the mini-batch operators. In particular, from Theorem~\ref{thm: comparison} we evince that smaller batch-sizes may lead to a better convergence rate than bigger batch-sizes. As underlined in Corollary~\ref{cr: comparison}, these theoretical results are in line with the literature on the convergence of the Bellman and the Gauss-Seidel operators. 
\begin{theorem}
\label{thm: comparison}
Consider a function $J:\mathcal{S}\rightarrow \mathbb{R}$ such that
$J(i) \leq (TJ)(i) \leq J^*(i)$ $\forall i\in\mathcal{S}$. Then, for any pair of integers $m$ and $m'$ such that $1\leq m'\leq m\leq n$, the following holds
\begin{equation}
\label{eq:inequality_th_6}
\begin{split}
&(B_{m}^kJ)(i) \leq (B_{m'}^kJ)(i) \leq J^*(i)\,,\quad i=1,\dots,n\,,\\
& \phantom{(B_{m}^kJ)(i) \leq (B_{m'}^kJ)(i) \leq J^*(i)\,,\quad } k=1,2,\dots\,.
\end{split}
\end{equation}

\begin{proof}
Similarly to Lemma~\ref{lmm: monotonicity} and Lemma~\ref{lmm: shift invariance}, we use induction twice, first on the state counter $i$ and then on the iteration counter $k$.
We set $k=1$ and apply induction over the state counter $i$; to simplify notation we drop the dependency on $k$.
As base-case we consider $i=1$. Since $\mathcal{M}_{m}(1)=\mathcal{M}_{m'}(1)=\{\emptyset\}$, then $(B_{m}J)(1)=(B_{m'}J)(1)=(BJ)(1)$. We now apply the $B$ operator on both sides of inequality $J(1)\leq J^*(1)$. The inequality $(BJ)(1) \leq J^*(1)$ follows directly from the monotonicity of the $B$ operator and the fact that $J^*$ is the unique fixed point of $B$, \textit{i.e.}, $(BJ^*)(1)=J^*(1)$.
Then we consider an arbitrary state index $i>1$ and  assume that $J(j) \leq (B_{m}J)(j) \leq (B_{m'}J)(j) \leq J^*(j)$ for $j=1,\dots,i-1$.
Since, without loss of generality, all states are assumed to be processed in ascending order and since all states in $\mathcal{M}(i)$ have already been processed by definition, then $\mathcal{M}(i) \subseteq \left\{1,2,\dots,i-1 \right\}$. Therefore, the induction assumption implies
\begin{equation*}
J(j) \leq (B_{m}J)(j) \leq (B_{m'}J)(j) \leq J^*(j)\,,\,\, \forall\,  j \in \mathcal{M}_{m}(i)\,,
\end{equation*}
from which it follows that
\begin{equation}
\label{eq: bound B_m}
0 \leq (B_{m}J)(j) - J(j) \leq (B_{m'}J)(j)- J(j)\,,\,\, \forall\,  j \in \mathcal{M}_{m}(i)\,.
\end{equation}
We now apply the $B_{m}$ operator to $J$ for state $i$ and using Equation~\eqref{eq: bound B_m} together with the fact that $\mathcal{M}_{m}(i) \subseteq \mathcal{M}_{m'}(i)$ by definition, we obtain
\begin{equation}
\label{eq: ineq comparative}
\begin{split}
&(B_m J)(i) \\&= \min_{u \in \mathcal{U}(i)} \!\!\left[  g(i,u) + \alpha\!\!\!\mkern-25mu\sum_{j \in \mathcal{S}\setminus \mathcal{M}_m(i) }\mkern-20mu\!\!\!\! p_{ij}(u) J(j) + \alpha\!\!\!\mkern-15mu \sum_{j \in \mathcal{M}_m(i)}\mkern-15mu\!\!\! p_{ij}(u) (B_mJ)(j)\right]\\
&= \!\!\min_{u \in \mathcal{U}(i)} \!\!\Bigg[  g(i,u) + \alpha\mkern-4mu \sum_{j \in \mathcal{S}}\mkern-2mu p_{ij}(u) J(j) + \alpha\!\!\mkern-15mu \sum_{j \in \mathcal{M}_m(i)}\mkern-15mu\!\! p_{ij}(u) (B_m\!-I)J(j)\!\Bigg] \\
&\leq \!\!\min_{u \in \mathcal{U}(i)} \!\!\Bigg[  g(i,u) + \alpha\!\mkern-4mu \sum_{j \in \mathcal{S}}\mkern-2mu p_{ij}(u) J(j) + \alpha\!\!\!\mkern-15mu \sum_{j \in \mathcal{M}_{m}(i)}\mkern-15mu\!\! p_{ij}(u)(B_{m'}\!-I)J(j)\!\Bigg]\\
&\leq\!\! \min_{u \in \mathcal{U}(i)} \!\!\Bigg[  g(i,u) + \alpha \mkern-4mu \sum_{j \in \mathcal{S}}\mkern-2mu p_{ij}(u) J(j) + \!\alpha\!\!\!\!\mkern-15mu \sum_{j \in \mathcal{M}_{m'}(i)}\mkern-15mu\!\!\! p_{ij}(u) (B_{m'}\!-I)J(j)\!\Bigg]\\ 
&=\min_{u \in \mathcal{U}(i)} \!\left[  g(i,u) + \alpha\!\!\!\! \mkern-25mu \sum_{j \in \mathcal{S}\setminus \mathcal{M}_{m'}(i) }\mkern-25mu\!\!\! p_{ij}(u) J(j) + \alpha \mkern-15mu\!\!\!\! \sum_{j \in \mathcal{M}_{m'}(i)}\mkern-15mu \!\!\!p_{ij}(u) (B_{m'}J)(j)\!\right]\\
&= (B_{m'}J)(i)\,,
\end{split}
\end{equation}
which is a valid induction step over $i$.

We now apply induction over the iteration counter $k$. The base-case is given by Equation~\eqref{eq: ineq comparative}. We assume that the following inequalities hold for an arbitrary integer $k\geq 1$ 
\begin{equation}
\label{eq: induction assumption k comparative}
(B_m^k J)(i)\leq (B_{m'}^kJ)(i)\,,\quad  i=1,\dots,n\,.
\end{equation}
Analogously to Lemma~\ref{lmm: monotonicity} and Lemma~\ref{lmm: shift invariance}, we introduce $\tilde{J}_m$ and $\tilde{J}_{m'}$ where $\tilde{J}_m = B_{m}^k J$ and $\tilde{J}_{m'} = B_{m'}^k J$. By applying the $B_m$ operator on both sides of Inequality~\eqref{eq: induction assumption k comparative} and making use of Inequality~\eqref{eq: ineq comparative} and the monotonicity property of $B$ (Lemma~\ref{lmm: monotonicity}), we obtain
\begin{equation*}
\begin{aligned}
(B_m^{k+1}J)(i) &= (B_m(B_{m}^kJ ))(i)\\
&= (B_m\tilde{J}_{m})(i)\\
&\leq(B_m \tilde{J}_{m'})(i)\\
&\leq(B_{m'} \tilde{J}_{m'})(i)\\
&=(B_{m'}(B_{m'}^kJ ))(i)\\
&= (B_{m'}^{k+1}J)(i)\,,
\end{aligned}
\end{equation*}
which is a valid induction step over $k$.
We can therefore conclude that $(B^{k}_{m}J)(i) \leq (B^{k}_{m'}J)(i)\leq J^*(i)$ for $i=1,\dots,n$ and $k=1,2,\dots$, where the latter inequality follows directly from the application of Lemma~\ref{lmm: monotonicity} and Proposition~\ref{prp: convergence}.
\end{proof}
\end{theorem}
\begin{corollary}
\label{cr: comparison}
Consider a function $J:\mathcal{S}\rightarrow \mathbb{R}$ such that
$J(i) \leq (TJ)(i) \leq J^*(i)$ $\forall i\in\mathcal{S}$. Then, for any integer $m$ such that $1\leq m \leq n$, the following holds
\begin{equation*}
\begin{split}
&(T^k J)(i) \leq (B_{m}^kJ)(i) \leq (F^kJ)(i)  \leq J^*(i)\,,\,\, i=1,\dots,n\,,\\
&\phantom{(T^k J)(i) \leq (B_{m}^kJ)(i) \leq (F^kJ)(i)  \leq J^*(i)\,,\,\,} k=1,2,\dots\,. 
\end{split}
\end{equation*}
\begin{proof}
The $T$ and $F$ operators are instances of the $B$ operator for specific values of the batch-size, \textit{i.e.}, $m=n$ for the $T$ operator and $m=1$ for the $F$ operator. The final result follows directly from the application of Theorem~\ref{thm: comparison} to these scenarios.
\end{proof}
\end{corollary}
\section{Numerical Evaluation}\label{sec:experiments}
This section is dedicated to numerically evaluating the mini-batch operators and the impact of different batch-sizes on the convergence of the main DP methods in practice when a multi-core system such as a GPU is deployed. 
\subsection{Set-Up and Environments}\label{sec:setup-and-environments}
For the accelerated tensor computations via graphic processing units (GPU), Pytorch~\cite{pytorch19} via the Python interface is deployed. All benchmarks are run using an NVIDIA GeForce GTX 1650 Ti and an Intel(R) Core(TM) i7-10750H CPU @ 2.60GHz. The code is publicly available at \url{https://gitlab.ethz.ch/gmatilde/mb-operator.git} in the form of a Python package. The performance of MB-VI and MB-MPI are studied on the following three OpenAI gym environments~\cite{openaigym16}, with a discount factor of 0.95 and $J(i)$ initialized to zero for all $i$. Before every evaluation of the operators, the states are shuffled and then selected and processed in blocks of $m$. 
\newline
\noindent{\textbf{FrozenLake:}}
The FrozenLake environment~\cite{openaigym16} consists of a $8\times 8$ grid-world with two types of tiles: normal tiles, whose associated cost is $1$, and high-penalty hole tiles, which have an associated cost of $10^3$. The goal is to go from the starting tile to the goal tile with the minimum cost. In each normal tile, the agent gets to choose among four possible directions, which will then lead the agent with a given probability to a next tile. If the agent ends up in a high-penalty tile, it will stay there with probability 1, consequently accumulating an infinite cost. The state and action spaces have dimension $64$ and $4$, respectively.   

\noindent{\textbf{Taxi:}}
The Taxi environment~\cite{dietterich00} consists of a $5\times 5$ grid with walls, where a taxi navigates picking-up and dropping-off passengers at given locations. The state is determined by the position of the taxi in the grid as well as the coordinates of the pick-up and drop-off locations of the passenger. A negative cost of -20 is given everytime a passenger is picked-up or dropped-off in the right location while illegal drop-off and pick-up actions have a cost of 10 and all other actions have a cost of 1. The transitions are deterministic and dictated by the geometry of the grid and the location of the walls. In addition, if the taxi plays a non-admissible action (one leading into a wall), it will remain at the same location. The state and action spaces have dimension 500 and 6, respectively (for a more detailed description see~\cite{dietterich00}). 

\noindent{\textbf{2D-Maze:}}
The 2D-Maze environment consists of a $N \times N$ maze with obstacles and the goal of the agent is to navigate to the terminal state with minimum cost. The transitions are not deterministic but every state-action-next state admissible triplet have an associated non-null probability. All state-action pairs have unitary cost, except for those involving the terminal state which have zero cost. For the benchmarks we use $N=80$ and $N=100$ which lead to state spaces with dimensions 6166 and 9706, respectively; the remaining grid is occupied by the maze walls. In both scenarios, the action space has dimension 4.  

\subsection{Benchmarks}\label{sec:benchmarks}

\begin{figure}[ht]
\begin{subfigure}{.2557\textwidth}
  \centering
  \includegraphics[width=0.98\linewidth]{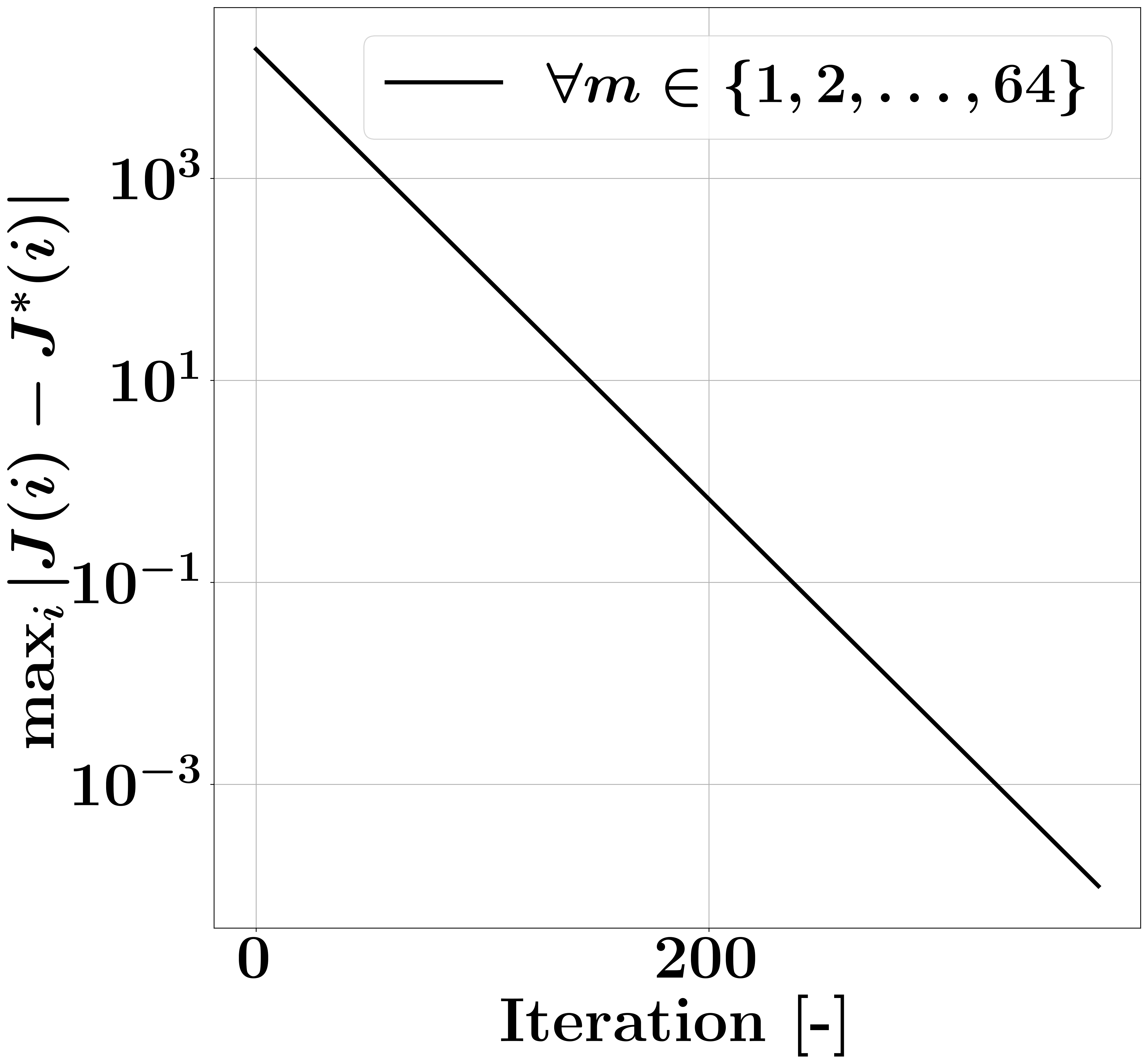}  
  \caption{Error vs iterations.}
  \label{fig:frozenlake_iterations}
\end{subfigure}\hspace{-0.17cm}
\begin{subfigure}{.242\textwidth}
  \centering
  \includegraphics[width=0.98\linewidth]{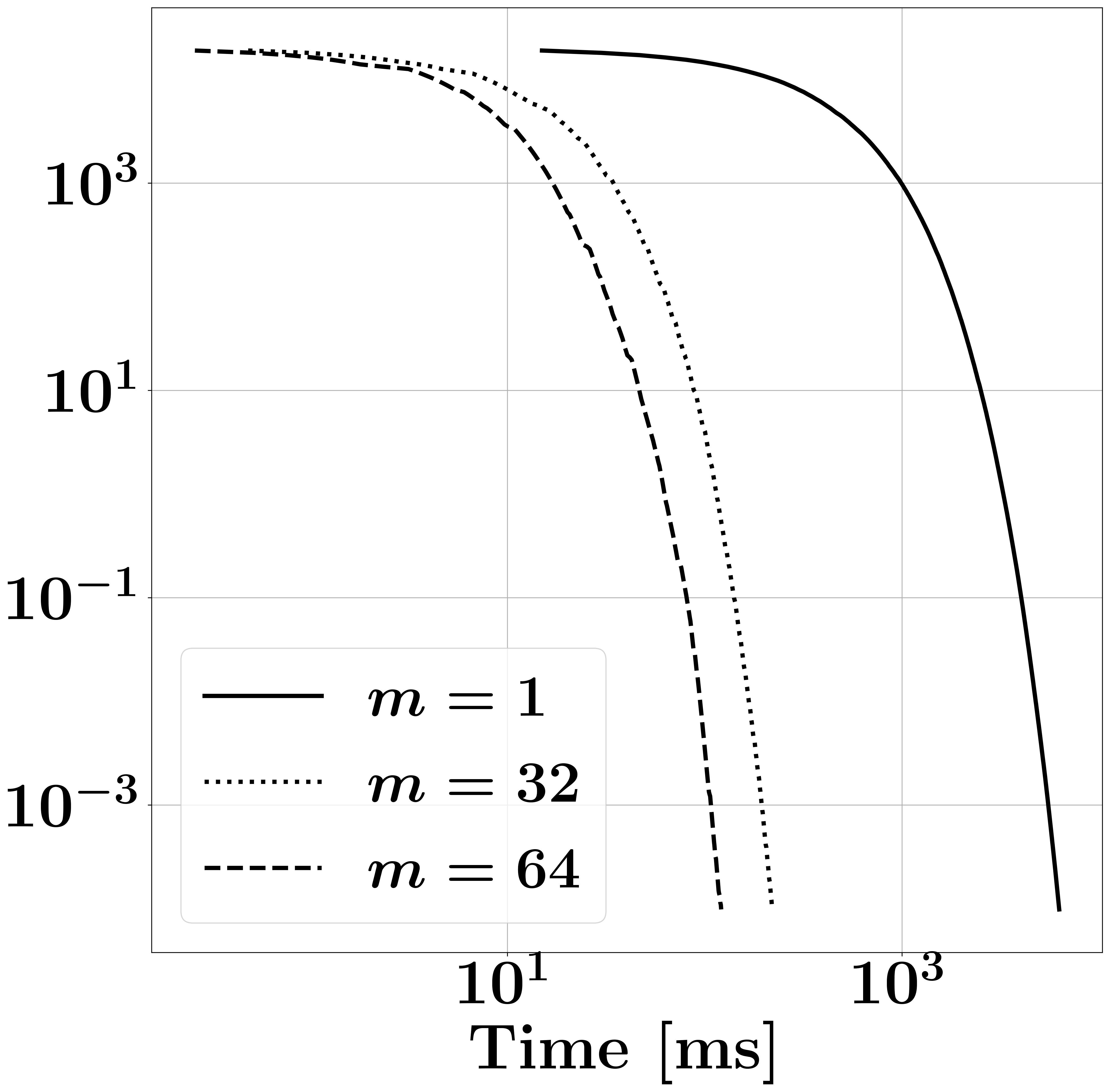}  
  \caption{Error vs computation time.}
  \label{fig:frozenlake_time}
\end{subfigure}
\caption{FrozenLake environment. We compare the convergence in terms of iterations and GPU time of VI ($m=64$), MB-VI ($m=32$) and GS-VI ($m=1$).}
\label{fig:frozenlake}
\end{figure}

\begin{figure}[ht]
\begin{subfigure}{.2557\textwidth}
  \centering
  \includegraphics[width=0.98\linewidth]{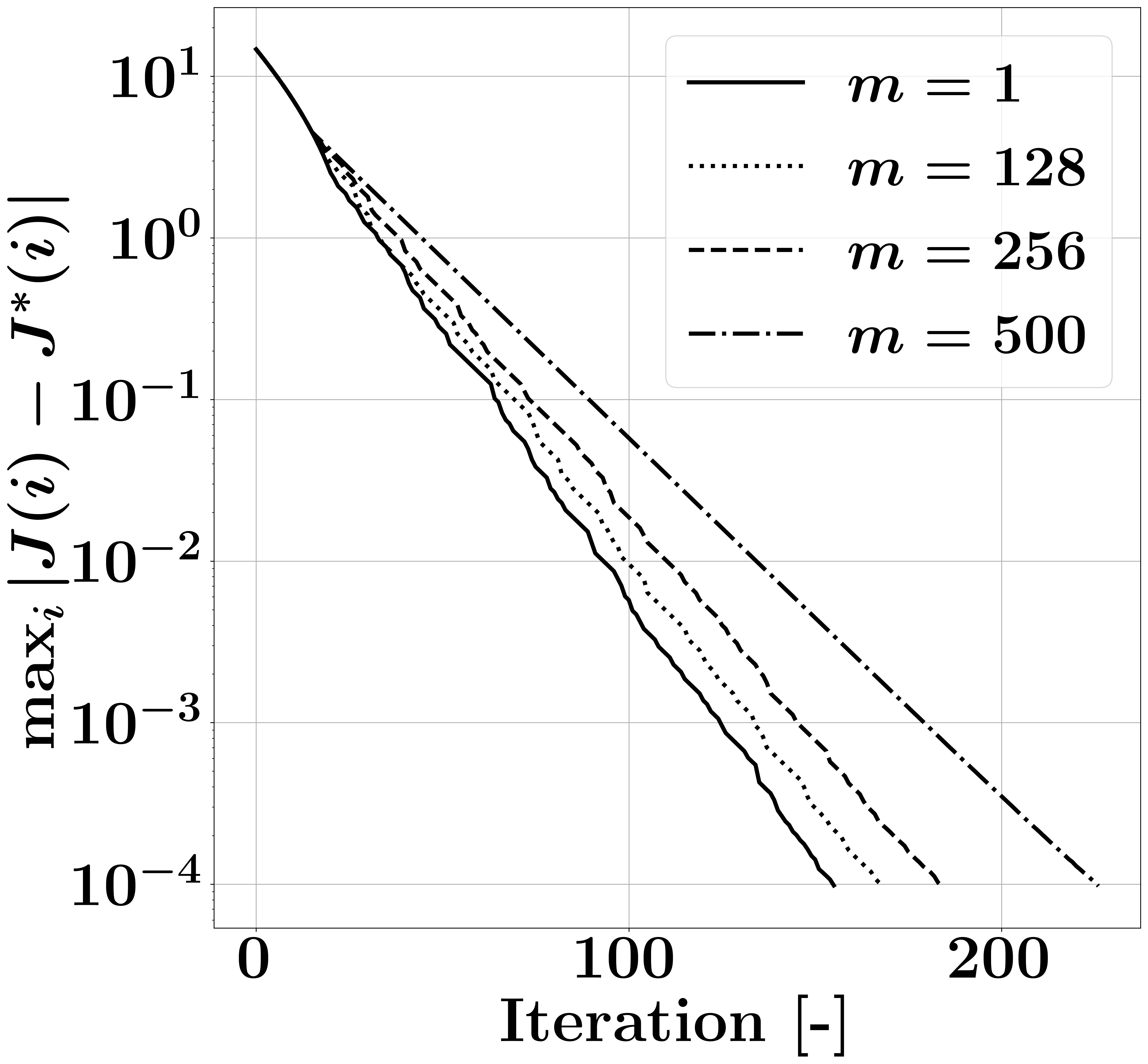}  
  \caption{Error vs iterations.}
  \label{fig:taxi_iterations}
\end{subfigure}\hspace{-0.17cm}
\begin{subfigure}{.242\textwidth}
  \centering
  \includegraphics[width=0.98\linewidth]{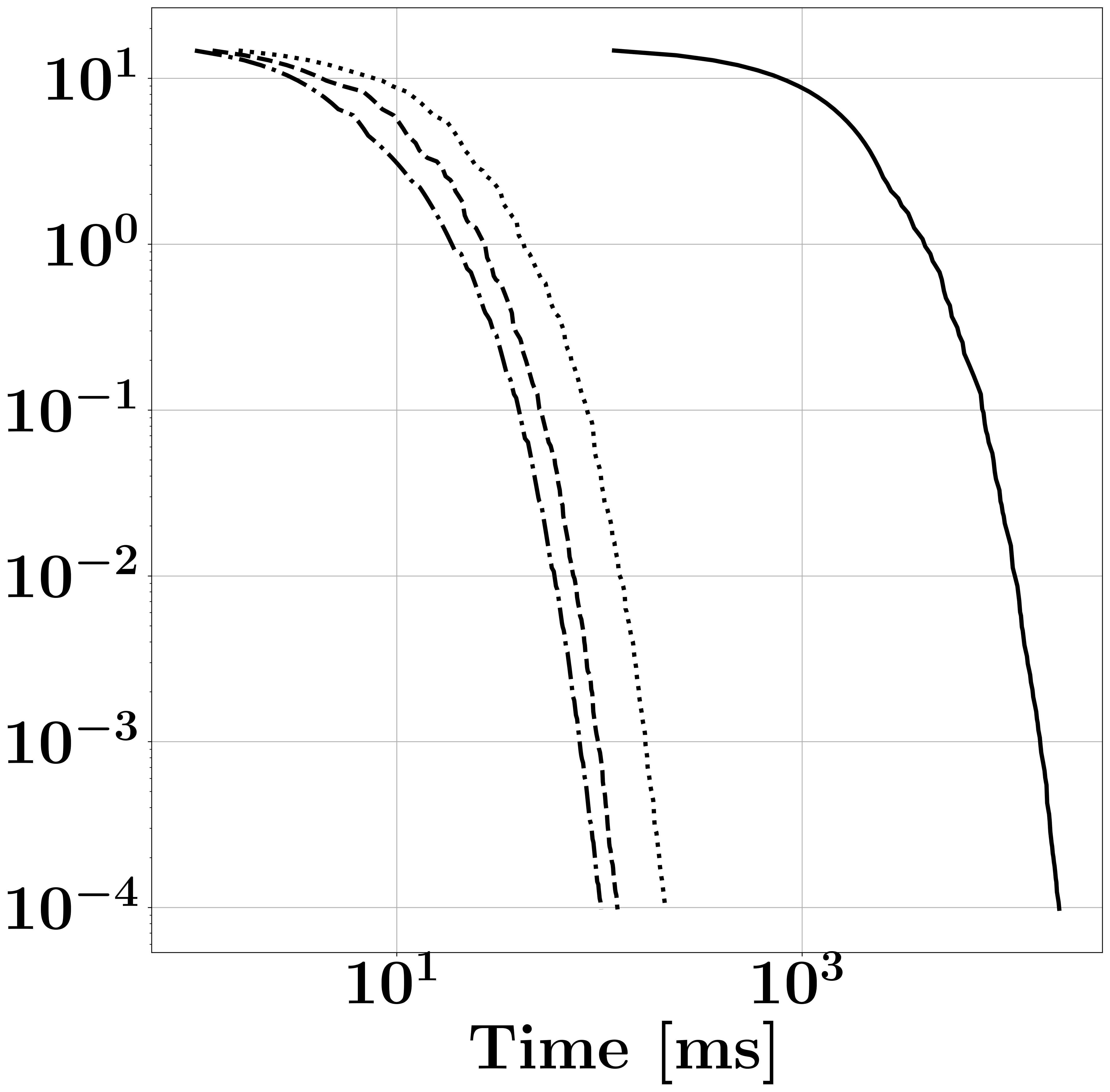}  
  \caption{Error vs computation time.}
  \label{fig:taxi_time}
\end{subfigure}
\caption{Taxi environment. We compare the convergence in terms of iterations and GPU time of VI ($m=500$), MB-VI ($m=128$, $m=256$) and GS-VI ($m=1$).}
\label{fig:taxi}
\end{figure}

\begin{figure}[ht]
\begin{subfigure}{.2557\textwidth}
  \centering
  \includegraphics[width=0.98\linewidth]{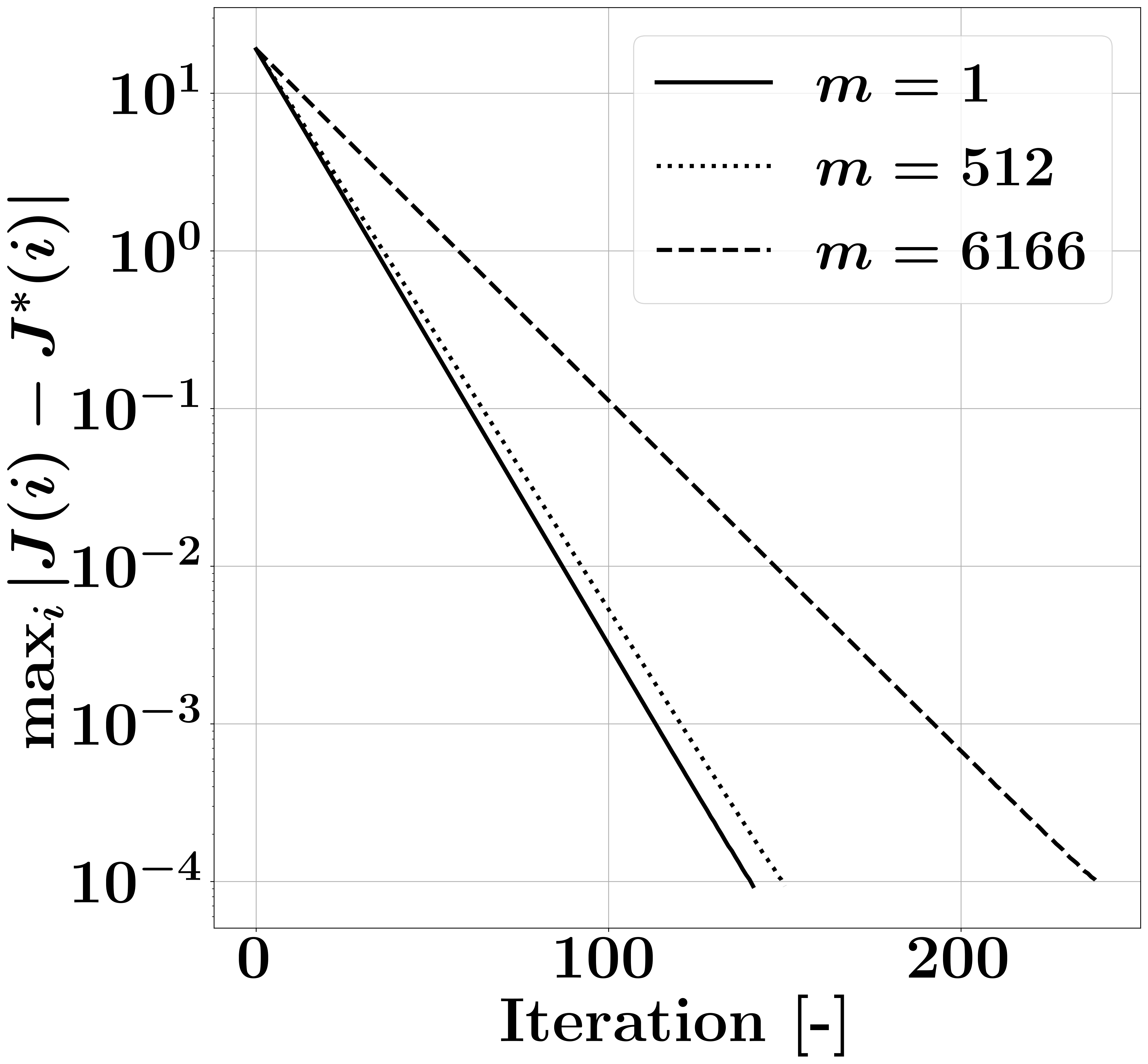}  
  \caption{Error vs iterations.}
  \label{fig:maze80_iterations}
\end{subfigure}\hspace{-0.17cm}
\begin{subfigure}{.242\textwidth}
  \centering
  \includegraphics[width=0.98\linewidth]{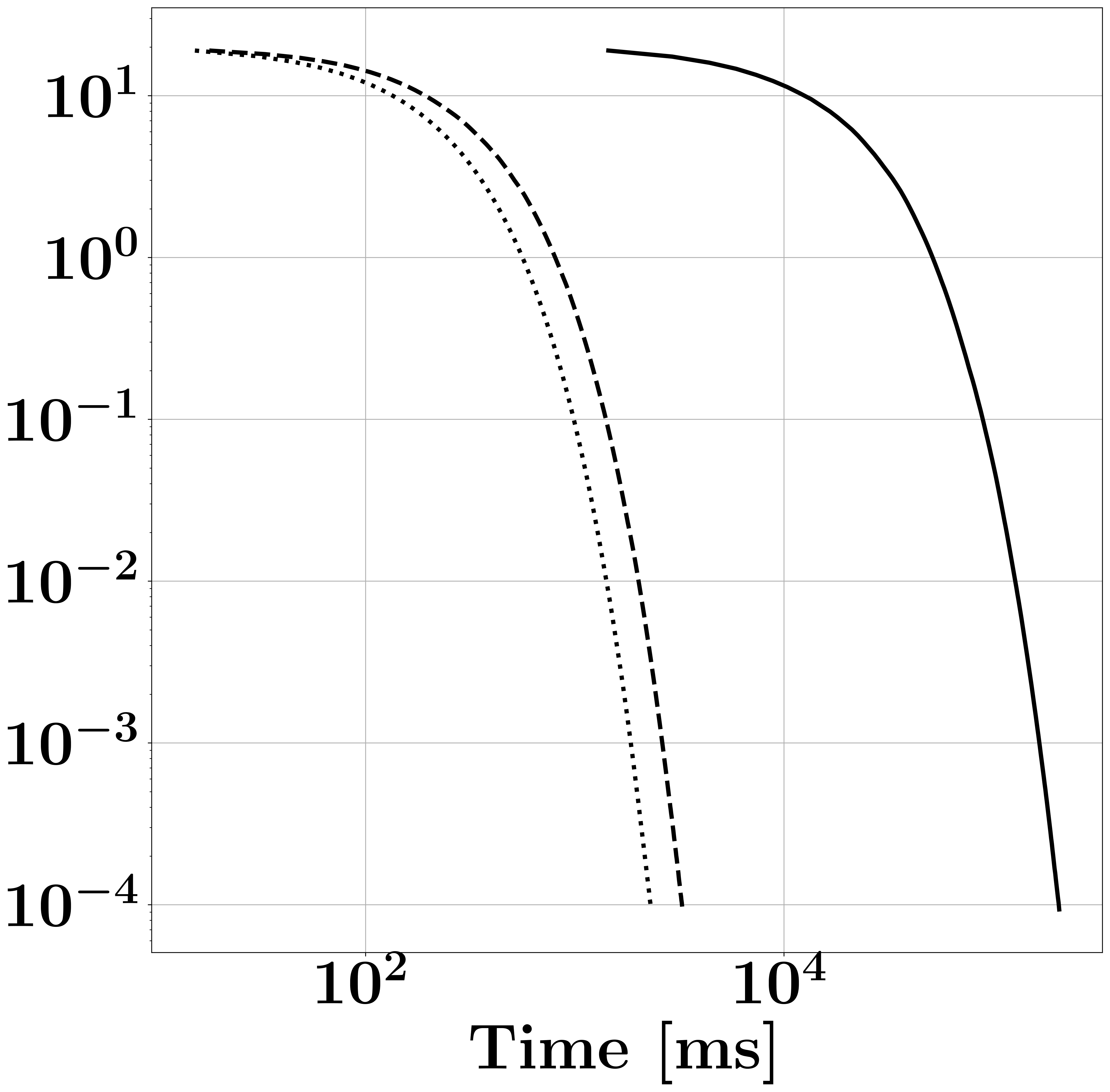}  
  \caption{Error vs computation time.}
  \label{fig:maze80_time}
\end{subfigure}
\caption{2D-Maze environment, $N=80$. We compare the convergence in terms of iterations and GPU time of VI ($m=6166$), MB-VI ($m=512$) and GS-VI ($m=1$).}
\label{fig:maze80}
\end{figure}

\begin{figure}[ht]
\hspace{-0.2cm}\begin{subfigure}{.2557\textwidth}
  \centering
  \includegraphics[width=\linewidth]{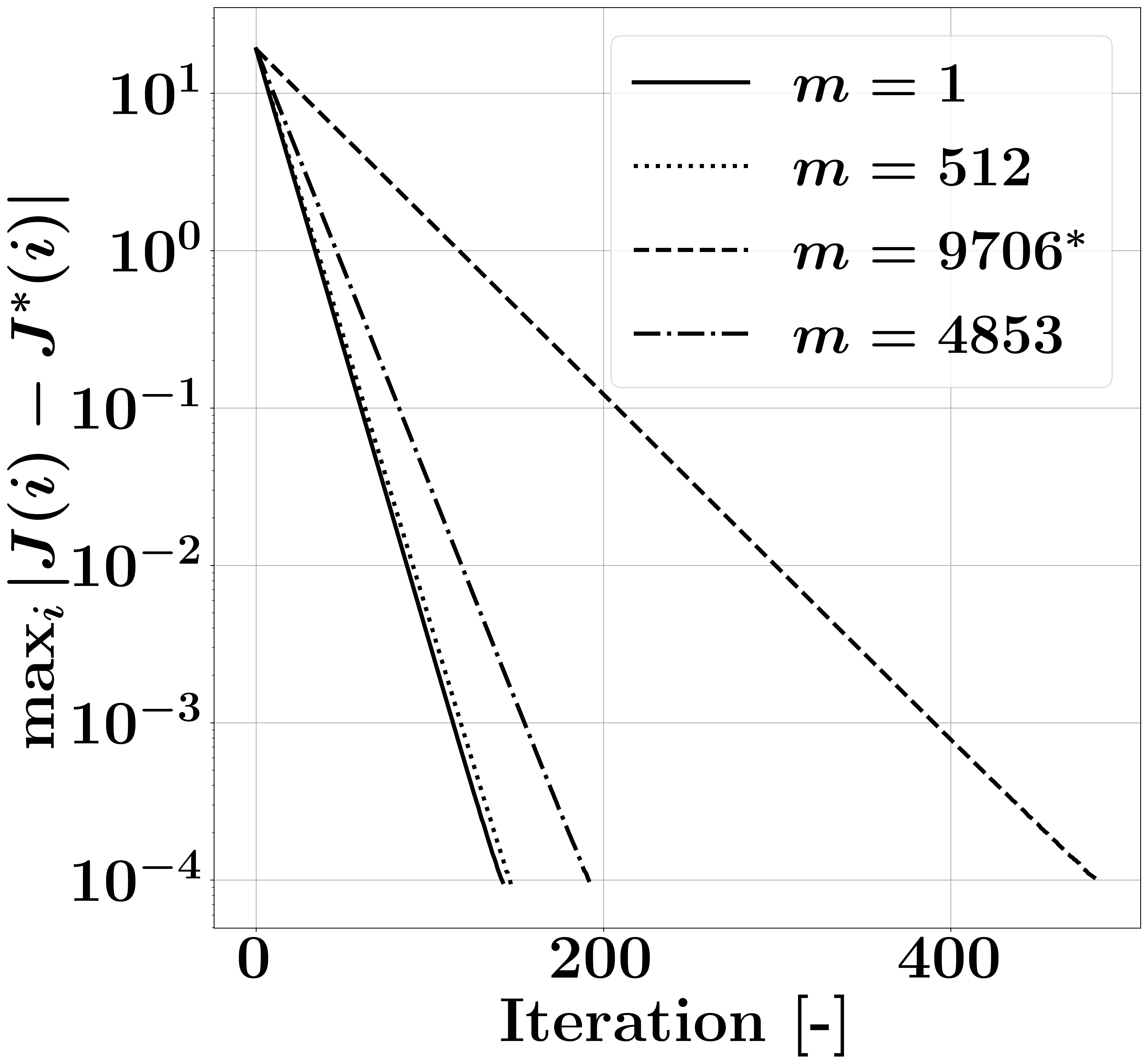}  
  \caption{Error vs iterations.}
  \label{fig:maze100_iterations}
\end{subfigure}\hspace{-0.11cm}
\begin{subfigure}{.242\textwidth}
  \centering
  \includegraphics[width=\linewidth]{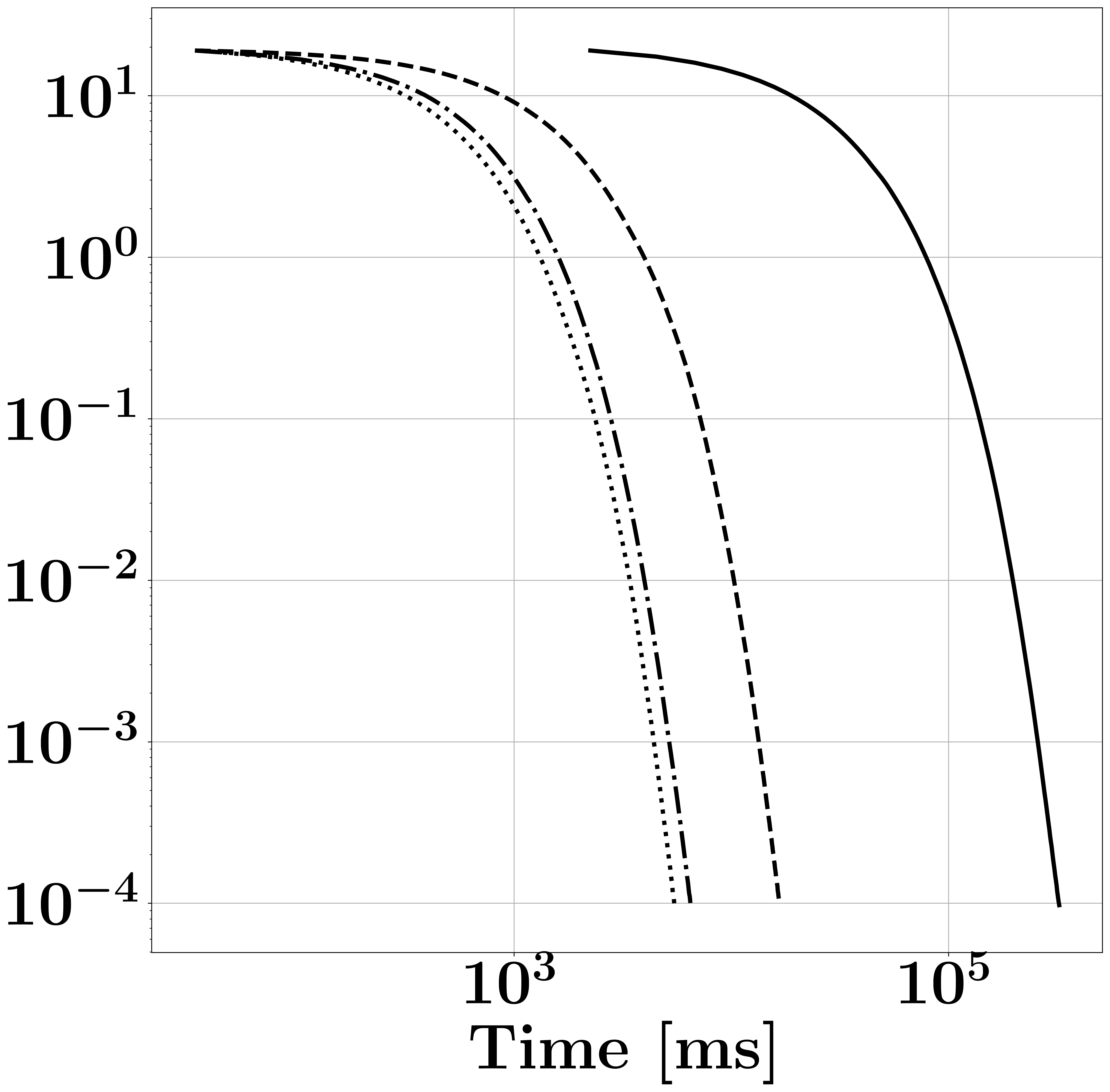}  
  \caption{Error vs computation time.}
  \label{fig:maze100_time}
\end{subfigure}
\caption{2D-Maze environment, $N=100$.  We compare the convergence in terms of iterations and GPU time of VI ($m=9706^*$), MB-VI ($m=512$, $m=4853$) and GS-VI ($m=1$). The asterisk indicates that the states are processed in batches of $4853$ because of the memory limitations of the GPU.}
\label{fig:maze100}
\end{figure}

For all the environments, we run MB-VI for different values of the batch-size $m$ and plot the infinity norm of the error versus number of iterations and GPU time. The plots in Figures~\ref{fig:frozenlake_iterations},~\ref{fig:taxi_iterations},~\ref{fig:maze80_iterations} and~\ref{fig:maze100_iterations} support the theoretical results of Theorem~\ref{thm: comparison} and Corollary~\ref{cr: comparison}. As we can see from these plots, the smaller is the batch-size and the faster the method converges in terms of number of iterations; this is clearly visible in Figures~\ref{fig:taxi_iterations} and~\ref{fig:maze80_iterations}, where the VI method ($m=500$ for the Taxi environment and $m=6166$ for the 2D-Maze environment with $N=80$) takes 71 and 98 iterations more than GS-VI ($m=1$) to converge to $10^{-4}$ of the optimal solution, respectively. Of course, as also stated in Theorem~\ref{thm: comparison}, there are settings where Inequality~\eqref{eq:inequality_th_6} is tight and therefore holds as equality. This is the case for the plot in Figure~\ref{fig:frozenlake_iterations}, where the convergence rate is not affected by the batch-size. Consequently, if we neglect the parallelization aspect, according to both the theoretical and the empirical results, the GS-VI method is the most effective as it always guarantees the fastest convergence in terms of number of iterations.

The convergence in terms of number of iterations is of course not the only aspect that should be taken into consideration when designing a numerical method. Modern hardware architectures, such as GPUs, offer the possibility of dramatically speeding up computation via massive parallelization~\cite{volkov08, schubiger20}; the parallelization capability of a method is relevant in determining its efficiency in terms of computation time~\cite{bertsekas97}. In this perspective, the iterations of the GS-VI method are not parallelizable. This ineherent sequentiality prevents one from the exploitation of parallel architectures and might result in longer overall computation times, despite the faster convergence in terms of number of iterations. On the other hand, the computation of the mini-batch operator for all the states within a batch is fully parallelizable; the extreme instance of this is the Bellman operator, whose evaluation is fully parallelizable. 
Because of its definition, the mini-batch operator offers flexibility in terms of the degree of parallelization which might lead to a speedup with respect to both the VI and GS-VI methods. All the reported speedup are computed for an accuracy of $10^{-4}$.

In the small-size scenarios in Figures~\ref{fig:frozenlake_iterations} and~\ref{fig:taxi_iterations}, given the hardware at hand, the best performing method is VI, achieving a speedup over GS-VI of $\times 51.71$ and $\times 182.23$, respectively, while with a smaller batch-size ($m=32$ for the FrozenLake environemnt and $m=256$ for the Taxi environment) we achieve a speedup of $\times 28.49$ and $\times 151.21$ over GS-VI, respectively. For the 2D-Maze environment with $N=80$ the best performance in terms of time is achieved when $m=512$. For this value of batch-size, we obtain a speedup of $\times 1.41$ over VI and of $\times 89.28$ over GS-VI. For the 2D-Maze environment with $N=100$, the full batch-size $m=9706$ can not be processed in parallel because of the memory limitations of the GPU. We therefore compute the Bellman operator in each iteration splitting the state space into two subsets of size $4853$ instead of processing all the states in one-shot and using the old values for the updates (this is indicated in the plot legend via an asterisk). We also run the MB-VI for $m=4853$ to see the advantages of MB-VI versus VI when both the methods enjoy the same parallelization degree. In general terms, when dealing with memory limitations, the mini-batch operator comes handy as we can reduce the batch-size in order to respect the memory limitations. This allows us to take advantage of the parallelization and, as shown in Figure~\ref{fig:maze100_iterations}, of faster convergence. With this environment the best time-performance is achieved with $m=512$, with a speedup of $\times 3.04$ over VI and $\times 58.86$ over GS-VI, while with $m=4853$ we have a speedup of $\times 2.57$ over VI and of $\times 49.73$ over GS-VI.

\begin{figure}[ht]
\hspace{-0.2cm}\begin{subfigure}{.2557\textwidth}
  \centering
  \includegraphics[width=\linewidth]{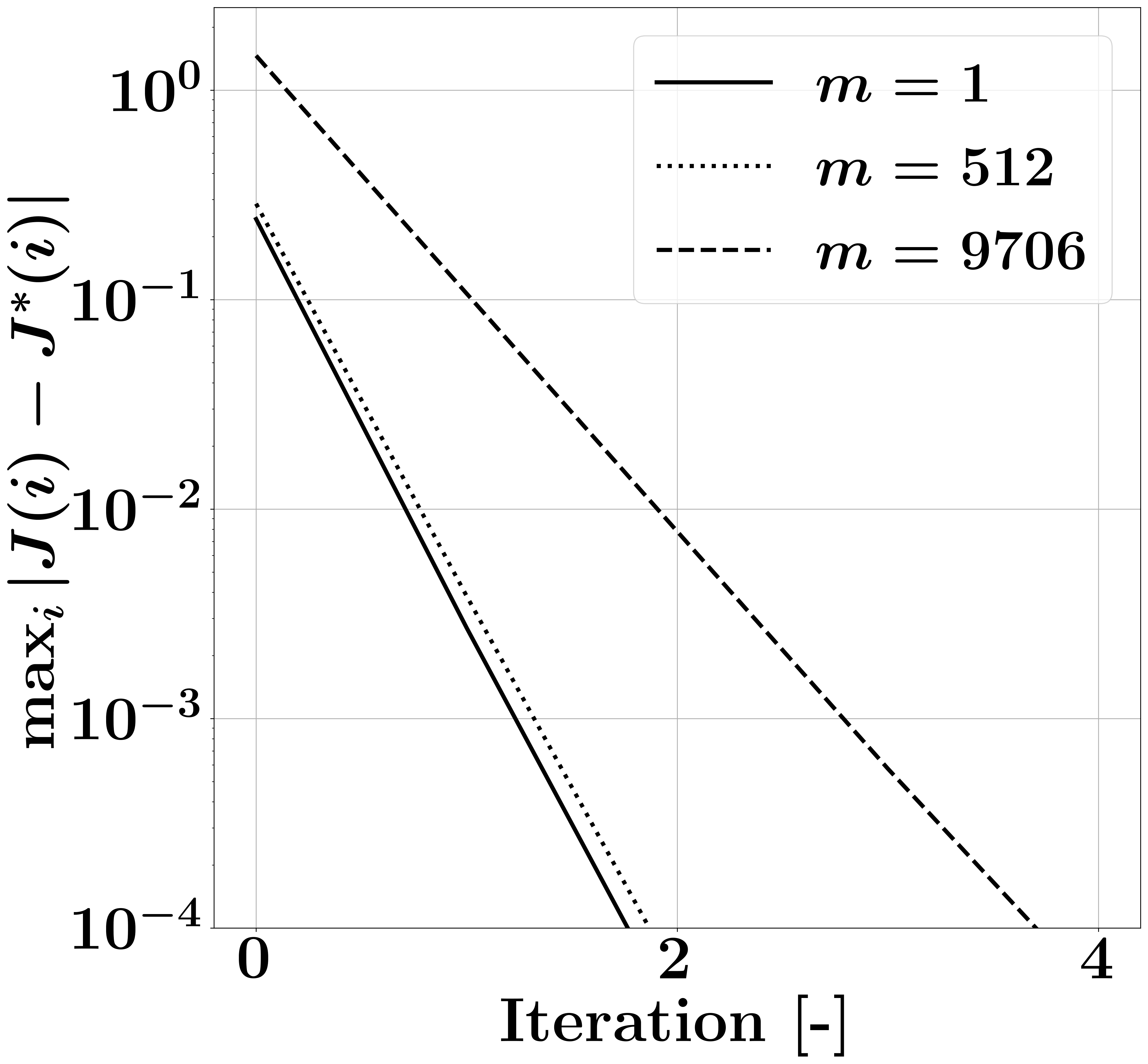}  
  \caption{Error vs iterations.}
  \label{fig:mazeMPI_iterations}
\end{subfigure}
\begin{subfigure}{.242\textwidth}
  \centering
  \includegraphics[width=\linewidth]{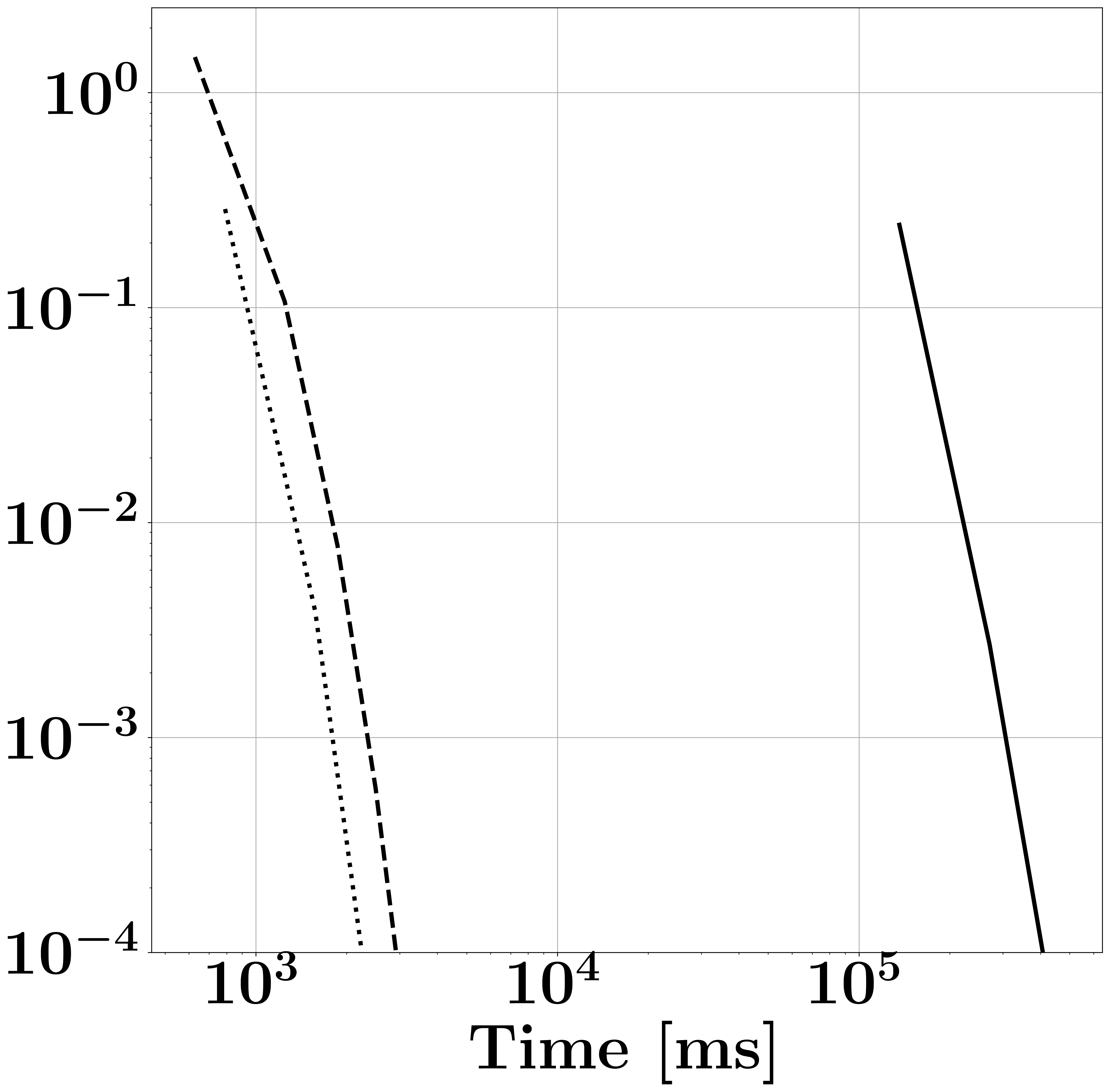}  
  \caption{Error vs computation time.}
  \label{fig:mazeMPI_time}
\end{subfigure}
\caption{2D-Maze environment, $N=100$.  We compare the convergence in terms of iterations and GPU time of MPI ($m=9706$), MB-MPI ($m=512$) and GS-MPI ($m=1$).}
\label{fig:mazeMPI}
\end{figure}

Similar considerations regarding the performance for different values of the batch-size also hold for the MB-MPI method with warm-start. As depicted in Figures~\ref{fig:mazeMPI_iterations} and~\ref{fig:mazeMPI_time}, for the 2D-Maze environment with $N=100$ and $K=50$ (see Algorithm~\ref{alg: MPI}) the best convergence in terms of number of iterations is achieved when $m=1$ while bigger batch-sizes require increasingly more iterations to achieve convergence. On the other hand, the best trade-off between convergence rate and parallelization speedup is achieved when $m=512$. For this batch-size the MB-MPI method converges $\times 1.31$ faster than MPI ($m=9706$) and $\times 194.52$ faster than GS-MPI ($m=1$).

\section{Conclusions \& Future Work}
We have proposed a new operator for dynamic programming based on the concept of mini-batch. The new operator generalizes the existing and well-enstablished Bellman and Gauss-Seidel operators for dynamic programming. When deployed in DP methods such as value iteration and modified policy iteration, our operator offers greater flexibility as the batch-size can be adjusted based on the hardware at hand as well as the specific problem structure. This allows the user to select the batch-size that realizes the best trade-off between convergence rate and parallelization speedup. We have proved that the new operator retains the key properties of shift-invariance, monotonicity and contractivity and we have characterized its fixed points. We have also characterized the convergence of the DP method based on the batch-size, showing that the deployment of bigger batch-sizes might result in slower rates of convergence. Finally, our empirical analysis demonstrates the competitive performance of MB-VI and MB-MPI across different benchmarks. 

For the time being, the randomization aspect of the new operator is not directly exploited. Even though the theoretical analysis is conducted under the assumption that the states are processed in ascending order, this of course does not have to be the case in practice. The definition of the mini-batch operator as well as the derived theoretical results also hold in case the states are uniformly sampled without replacement, as currently done for the benchmarks. There could, however, be advanteges in sampling the states with replacement and/or according to a non-uniform distribution. In this perspective, one can inverstigate the use of importance-sampling based techniques~\cite{zhao15} and $\epsilon$-greedy policy strategies~\cite{tokic10} for the selection of the states in the mini-batches. This would require further work on the theoretical side, but could potentially improve the convergence rates by exploiting the specific problem structure. 
Another interesting future line of research would be asynchronous variants of MB-VI and MB-MPI, which hold the promise of further speedup by avoiding synchronization~\cite{tsitsiklis86}.

\section*{Acknowledgment}
This work has been supported by the European Research Council (ERC) under the H2020
Advanced Grant no. 787845 (OCAL).


\begin{thebibliography}{00}
\bibitem{tsitsiklis86} J. N. Tsitsiklis, D. P. Bertsekas, and M. Athans, ``Distributed asynchronous deterministic stochastic gradient optimization algorithms,'' \emph{IEEE  Transaction on Automatic Control}, vol. 31, no. 9, pp 803-812, Sep. 1986.

\bibitem{rockafellar76} R. T. Rockafellar, ``Monotone operators and the proximal point algorithm,'' \emph{SIAM Journal on Control and Optimization}, vol. 14, no. 5, pp.877-898, Aug. 1976.

\bibitem{bellman52} R. Bellman, ``On the theory of dynamic programming,'' \emph{Proceedings of the National Academy of Sciences}, vol. 38, no. 8, pp. 716-719, Aug. 1952.

\bibitem{bobrow85} J. E. Bobrow, S. Dubowsky, and J. S. Gibson, ``Time-optimal
control of robotic manipulators along specified paths,'' \emph{International Journal of Robotics Research}, vol. 4, no. 3, pp. 3-17, Sep. 1985.

\bibitem{bertsimas98} D. Bertsimas, and A.W. Lo, ``Optimal control of execution
costs,'' \emph{Journal of Financial Markets}, vol. 1, no. 1, pp. 1-50, 1998.

\bibitem{beck13} A. Beck, and L. Tetruashvili, ``On the convergence of block coordinate descent type
methods,'' \emph{SIAM Journal on Optimization}, vol. 23, no. 4, pp. 2037-2060, Oct. 2013.

\bibitem{bottou18} L. Bottou, F. E. Curtis, and J. Nocedal, ``Optimization methods for large-scale machine learning,'' \emph{SIAM Review}, vol. 60, no. 2, pp. 223-311, May 2018.

\bibitem{schubiger20} M. Schubiger, G. Banjac, and J. Lygeros, ``{GPU} acceleration of {ADMM} for large-scale quadratic programming,'' \emph{Journal of Parallel and Distributed Computing}, vol. 144, pp. 55-67, Oct. 2020.

\bibitem{dietterich00} T. G. Dietterich, ``Hierarchical reinforcement learning with the MAXQ value function decomposition,'' \emph{Journal of Arti cial Intelligence Research}, vol. 13, pp. 227-303, Nov. 2000.

\bibitem{bellman57} R. Bellman, \emph{Dynamic Programming,} Princeton University Press, Princeton, New Jersey, USA, 1957.

\bibitem{sutton17} R. Sutton, \emph{Reinforcement Learning: An Introduction,} 2\textsuperscript{nd} ed., The MIT Press, Cambridge, USA, 2018.

\bibitem{bertsekas12} D. P. Bertsekas, \emph{Dynamic Programming and Optimal Control,} 4\textsuperscript{th} ed., vol. 2, Athena Scientific, Bellmont, USA, 2012.

\bibitem{bertsekas99} D. P. Bertsekas, \emph{Nonlinear Programming,} 2\textsuperscript{nd} ed., Athena Scientific, Bellmont, USA, 1999.


\bibitem{bertsekas97} D. P. Bertsekas, and J. N. Tsitsiklis, \emph{Parallel and Distributed Computation: Numerical Methods,} Athena Scientific, Bellmont, USA, 1997.

\bibitem{kelley95} C. T. Kelley, \emph{Iterative Methods for Linear and Nonlinear Equations,} Society  for  Industrial  and  Applied  Mathematics, Philadelphia, USA, 1995.


\bibitem{martinelli20} A. Martinelli, M. Gargiani, and J. Lygeros, ``Data-driven optimal control with a relaxed linear program,'' \emph{Automatica (to appear)}, 2021.

\bibitem{richtarik2011} P. Richt\'arik, and M. Takáč, ``Iteration complexity of randomized block-coordinate descent methods for minimizing a composite function,'' \emph{arXiv: 1107.2848}, 2011.

\bibitem{openaigym16} G. Brockman, V. Cheung, L. Pettersson, J. Schneider, J. Schulman, J. Tang, and W. Zaremba, ``OpenAI Gym,'' \emph{arXiv:1606.01540}, 2016.

\bibitem{pytorch19} A. Paszke, S. Gross, F. Massa, A. Lerer, J. Bradbury, G. Chanan, T. Killeen, Z. Lin, N. Gimelshein, L. Antiga, A. Desmaison, A. Kopf, E. Yang, Z. DeVito, M. Raison, A. Tejani, S. Chilamkurthy, B. Steiner, L. Fang, J. Bai, and S. Chintala, ``PyTorch: an imperative style, high-performance deep learning library,'' \emph{Advances in Neural Information Processing Systems 32}, 2019.

\bibitem{volkov08} V. Volkov, and J. W. Demmel, ``Benchmarking GPUs to tune dense linear algebra,'' \emph{SC '08: Proceedings of the 2008 ACM/IEEE Conference on Supercomputing}, 2008.

\bibitem{tokic10} M. Tokic, ``Adaptive $\epsilon$-greedy exploration in reinforcement learning based on value differences,'' \emph{KI 2010: Advances in Artificial Intelligence}, 2010.

\bibitem{zhao15} P. Zhao, and T. Zhang, ``Stochastic optimization with importance sampling for regularized loss minimization,'' \emph{Proceedings of the 32nd International Conference on Machine Learning}, 2015.

\end{thebibliography}
\end{document}